\documentclass[12pt]{amsart}

\title{On certain classes of modules over group algebras}
\author{Ioannis Emmanouil and Wei Ren}

\newtheorem{Lemma}{Lemma}[section]
\newtheorem{Proposition}[Lemma]{Proposition}
\newtheorem{Theorem}[Lemma]{Theorem}
\newtheorem{Corollary}[Lemma]{Corollary}

\begin{document}

\begin{abstract}
In this paper, we examine the relation between certain subclasses
of the classes of Gorenstein projective, Gorenstein flat and
Gorenstein injective modules over a group algebra, which consist
of the cofibrant, cofibrant-flat and fibrant modules respectively.
These subclasses have all structural properties that the Gorenstein
classes are known to have, regarding the existence of complete
cotorsion pairs and various approximations between them. Furthermore,
these subclasses have certain properties, which are conjecturally
anticipated to hold for the Gorenstein classes as well. Since the
Gorenstein classes are actually equal to the corresponding
subclasses for large families of groups, we thus obtain (indirectly)
some information about the Gorenstein classes over these groups.
\end{abstract}

\makeatletter
\@namedef{subjclassname@2020}{\textup{2020} Mathematics Subject Classification}
\makeatother

\thanks{}
\date{}


\maketitle
\tableofcontents

\addtocounter{section}{-1}
\section{Introduction}

\noindent
In homological algebra, flat modules play an important intermediary
role between projective and injective modules: Projective modules
are flat, whereas flat modules are those modules whose character
modules are injective. Gorenstein homological algebra is the relative
homological theory, which is based upon Gorenstein projective,
Gorenstein injective and Gorenstein flat modules; cf.\ \cite{EJ}.
The theory has developed rapidly during the past few decades and found
many applications in the representation theory of Artin algebras
and the theory of singularities. In contrast to their classical
homological algebra counterparts, it is not known though whether
(i) Gorenstein projective modules are Gorenstein flat and (ii)
Gorenstein flat modules are those modules whose character modules
are Gorenstein injective. These two questions represent fundamental
problems in Gorenstein homological algebra. Holm proved in \cite{Hol}
that both questions admit a positive answer if the ring is right
coherent and (for question (i)) has finite left finitistic dimension.
\v{S}aroch and \v{S}t$\!$'$\!$ov\'{i}\v{c}ek introduced in \cite{SS}
the class of projectively coresolved Gorenstein flat modules and
showed that these modules are Gorenstein projective as well. An
equivalent formulation of question (i) above is to ask whether any
Gorenstein projective module is projectively coresolved Gorenstein
flat. If this is true, then all flat Gorenstein projective modules
are actually projective; the latter property is studied in \cite{BCIE}.
We may obtain yet another equivalent formulation of (i) by modifying
(ii): All Gorenstein projective modules are Gorenstein flat if and
only if the character modules of Gorenstein projective modules are
Gorenstein injective.

Gorenstein homological algebra has also found applications in
cohomological group theory. The Gorenstein cohomological dimension
of a group $G$ over $\mathbb{Z}$ generalizes the ordinary cohomological
dimension of $G$ and coincides for certain classes of groups with the
geometric Bredon dimension, i.e.\ with the minimal dimension of a model
for the classifying space $\underline{E}G$ for proper actions of $G$;
see the Introduction of \cite{ET1} for some details on that relation.
The Gorenstein cohomological dimension of groups over $\mathbb{Z}$ is
proposed in \cite{BDT} to serve as an algebraic invariant, whose finiteness
characterizes the groups $G$ that admit a finite dimensional model for
$\underline{E}G$.

In order to test various questions regarding properties of the Gorenstein
module classes to the case where the ring is specialized to be the group
algebra of a group $G$ with coefficients in a commutative ring $k$, it is
useful to consider certain subclasses of the Gorenstein module classes,
that admit a simpler and more robust description. We denote by $B$ the
$kG$-module consisting of those $k$-valued functions on $G$ whose image
is a finite set. Benson defined in \cite{Ben} a $kG$-module $M$ to be
cofibrant if the diagonal $kG$-module $M \otimes_kB$ is projective. It
turns out that cofibrant modules are precisely the cokernels of those
acyclic complexes of projective $kG$-modules $P$, for which the complex
$P \otimes_kB$ is contractible. Analogously, we say that a $kG$-module
$M$ is cofibrant-flat if the diagonal $kG$-module $M \otimes_kB$ is flat;
these modules are precisely the cokernels of those acyclic complexes of
flat $kG$-modules $F$, for which the complex $F \otimes_kB$ is pure-acyclic.
In the same way, we say that a $kG$-module $M$ is fibrant if the diagonal
$kG$-module $\mbox{Hom}_k(B,M)$ is injective. Fibrant modules are precisely
the kernels of those acyclic complexes of injective $kG$-modules $I$, for
which the complex $\mbox{Hom}_k(B,I)$ is contractible. Elaborating on
a result by Cornick and Kropholler \cite{CK}, who showed that cofibrant
modules are Gorenstein projective, it follows that cofibrant-flat (resp.\
fibrant) modules are Gorenstein flat (resp.\ Gorenstein injective); cf.\
\cite{ET3}. In view of the very definition of cofibrant, cofibrant-flat
and fibrant modules, it is clear that (i) any cofibrant module is
cofibrant-flat and (ii) the cofibrant-flat modules are those $kG$-modules
whose character modules are fibrant. The class of cofibrant modules is the
left hand side of a complete hereditary cotorsion pair; cf.\ \cite{ER}.
Analogously, we show that the class of cofibrant-flat (resp.\ fibrant)
modules is the left (resp.\ right) hand side of a complete hereditary
cotorsion pair as well. Approximating cofibrant-flat modules by cofibrant
modules, we conclude that the class of cofibrant-flat modules is generated
(as an exact subcategory of the module category) by the class of cofibrant
and the class of flat modules.

The homological behaviour of the $kG$-module $B$ is manageable over
groups in Kropholler's class ${\scriptstyle{{\bf LH}}}\mathfrak{F}$
of hierarchically decomposable groups \cite{Kro}. We recall that
${\scriptstyle{{\bf H}}}\mathfrak{F}$ is the smallest class of groups
that contains all finite groups and is such that whenever a group $G$
admits a cellular action on a finite dimensional contractible CW-complex
with stabilizers in ${\scriptstyle{{\bf H}}}\mathfrak{F}$, then $G$ is
also contained in ${\scriptstyle{{\bf H}}}\mathfrak{F}$. The class
${\scriptstyle{{\bf LH}}}\mathfrak{F}$ contains those groups, all of
whose finitely generated subgroups are contained in
${\scriptstyle{{\bf H}}}\mathfrak{F}$; we note that
${\scriptstyle{{\bf LH}}}\mathfrak{F}$ is a big class of groups.
Over ${\scriptstyle{{\bf LH}}}\mathfrak{F}$-groups, the arguments
by Dembegioti and Talelli \cite{DT}, who showed that any Gorenstein
projective module is cofibrant, show that any Gorenstein flat (resp.\
Gorenstein injective) module is cofibrant-flat (resp.\ fibrant), at
least in the case where $k=\mathbb{Z}$ is the ring of integers; see
also \cite{ET3}. Then, the properties of cofibrant, cofibrant-flat
and fibrant modules provide us with corresponding results regarding
the Gorenstein module classes. To state such an example, we note that,
if $G$ is an ${\scriptstyle{{\bf LH}}}\mathfrak{F}$-group, then there
exists a $\mathbb{Z}G$-module $T$, such that the class of Gorenstein
flat $\mathbb{Z}G$-modules is precisely the kernel of the additive
functor $\mbox{Tor}_1^{\mathbb{Z}G}(T,\_\!\_)$; in particular, the
class of Gorenstein flat $\mathbb{Z}G$-modules is closed under pure
submodules and pure quotients. In fact, the class of Gorenstein flat
$\mathbb{Z}G$-modules is precisely the class of pure quotients of
Gorenstein projective $\mathbb{Z}G$-modules (if $G$ is an
${\scriptstyle{{\bf LH}}}\mathfrak{F}$-group).

Here is a description of the contents of the paper: Following the
preliminary Section 1, in Section 2, we introduce and study the
cofibrant-flat modules, which are related to Gorenstein flat
modules. In Section 3, we examine approximations of cofibrnat-flat
modules by cofibrant modules, that specialize the approximations of
Gorenstein flat modules by Gorenstein projective modules. Finally,
in Section 4, we examine the class of fibrant modules, relate these
modules to Gorenstein injective modules and study the duality between
cofibrant-flat and fibrant modules.

\vspace{0.1in}

\noindent
{\em Notations and terminology.}
We consider a commutative ring $k$ and a group $G$. Unless
otherwise specified, all modules are modules over the group
algebra $kG$. We denote by ${\tt Proj}(kG)$, ${\tt Inj}(kG)$,
${\tt Flat}(kG)$ and ${\tt Cotor}(kG)$ the classes of projective,
injective, flat and cotorsion $kG$-modules respectively. If
$M$ is any $kG$-module, then its character (Pontryagin dual)
module is the $kG$-module
$DM = \mbox{Hom}_{\mathbb{Z}}(M,\mathbb{Q}/\mathbb{Z})$ of
all additive maps from $M$ to $\mathbb{Q}/\mathbb{Z}$.

\section{Preliminaries}

\noindent
In this section, we record certain basic notions which are used
in the sequel. These concern modules over group algebras, cotorsion
pairs and model structures in module categories.

\vspace{0.1in}

\noindent
{\sc I.\ Modules over group algebras.}
Let $k$ be a commutative ring and $G$ a group. Using the diagonal
action of $G$, the tensor product $M \otimes_k N$ of two $kG$-modules
$M,N$ admits the structure of a $kG$-module, where
$g \cdot (x \otimes y) = gx \otimes gy \in M \otimes_k N$ for any
$g \in G$, $x \in M$ and $y \in N$. If the $kG$-module
$M$ is projective and $N$ is projective as a $k$-module, then the
diagonal $kG$-module $M \otimes_k N$ is projective as well; cf.\
\cite[Chapter III, Corollary 5.7]{Bro}. Analogously, if the
$kG$-module $M$ is flat and $N$ is flat as a $k$-module, then the
diagonal $kG$-module $M \otimes_k N$ is flat. It follows readily
that $\mbox{fd}_{kG}(M \otimes_kN) \leq \mbox{fd}_{kG}M$, whenever
$N$ is $k$-flat. We also note that for any two $kG$-modules $M,N$
the $k$-module $\mbox{Hom}_k(M,N)$ admits the structure of a
$kG$-module, with the group $G$ acting diagonally; here,
$(g \cdot f)(x) = gf(g^{-1}x)\in N$ for any $g \in G$,
$f \in \mbox{Hom}_k(M,N)$ and $x \in M$. The $G$-invariant part
of the $kG$-module $\mbox{Hom}_k(M,N)$ is precisely the $k$-module
$\mbox{Hom}_{kG}(M,N)$ of $kG$-linear maps from $M$ to $N$. If the
$kG$-module $M$ is projective and $N$ is injective as a $k$-module,
then the diagonal $kG$-module $\mbox{Hom}_k(M,N)$ is injective as well.
It follows that $\mbox{id}_{kG}\mbox{Hom}_k(M,N) \leq \mbox{pd}_{kG}M$,
whenever $N$ is $k$-injective. We note that for any $kG$-modules
$L,M,N$ the natural isomorphism of $k$-modules
\[ \mbox{Hom}_k(L \otimes_k M,N) \simeq
   \mbox{Hom}_k(L, \mbox{Hom}_k(M,N)) \]
is actually an isomorphism of $kG$-modules, where $G$ acts diagonally
on the tensor product and the three Hom-groups that are involved. Hence,
taking $G$-invariants, we conclude that the above isomorphism of
$kG$-modules restricts to an isomorphism of $k$-modules
\[ \mbox{Hom}_{kG}(L \otimes_k M,N) \simeq
   \mbox{Hom}_{kG} (L,\mbox{Hom}_k(M,N)) . \]
We conclude that the $kG$-module $\mbox{Hom}_k(M,N)$ is injective
if $M$ is $k$-flat and $N$ is injective.

Let $B$ be the $kG$-module consisting of all functions from $G$ to
$k$ that admit finitely many values, introduced in \cite{KT} in the
special case where $k = \mathbb{Z}$. The $kG$-module $B$ is $kH$-free
for any finite subgroup $H \subseteq G$. The constant functions induce
a $kG$-linear map $\iota : k \longrightarrow B$. In fact, $\iota$ is
a $k$-split monomorphism; a $k$-linear splitting for $\iota$ may be
obtained by evaluating functions at the identity element of $G$. We
also let $\overline{B} = \mbox{coker} \, \iota$. The analysis of the
homological properties of the $kG$-module $B$ is tractable over the
class ${\scriptstyle{{\bf LH}}}\mathfrak{F}$ of hierarchically
decomposable groups, that were defined  by Kropholler in \cite{Kro}.
All soluble groups, all groups of finite virtual cohomological
dimension, all one-relator groups and all automorphism groups of
Noetherian modules over commutative rings are
${\scriptstyle{{\bf LH}}}\mathfrak{F}$-groups. Moreover, the class
${\scriptstyle{{\bf LH}}}\mathfrak{F}$ is closed under
extensions, ascending unions, free products with amalgamation and
HNN-extensions. Another class of groups over which the homological
behaviour of $B$ can be controlled consists of the groups of type
$\Phi$, which were introduced in \cite{T}, in connection with the
existence of a finite dimensional model for the classifying space
for proper actions;  see also \cite{MS}. A group $G$ is of type
$\Phi$ over $k$ whenever the class of $kG$-modules of finite projective
dimension coincides with the class of those $kG$-modules that have
finite projective dimension over any finite subgroup of $G$. It follows
that $\mbox{pd}_{kG}B < \infty$, if $G$ is of type $\Phi$ over $k$.

\vspace{0.1in}

\noindent
{\sc II.\ Closure operations on group classes.}
The construction of ${\scriptstyle{{\bf LH}}}\mathfrak{F}$-groups
is a particular case for the class $\mathfrak{F}$ of finite groups
of a general operation that associates with any group class
$\mathfrak{C}$ the class ${\scriptstyle{{\bf LH}}}\mathfrak{C}$,
as defined by Kropholler \cite{Kro}. If $k$ is a fixed commutative
ring, then the definition by Talelli \cite{T} of groups of type
$\Phi$ over $k$ is a particular case of an operation on group classes
as well: If $\mathfrak{C}$ is a class of groups, we let
$\Phi \mathfrak{C}$ be the class consisting of those groups $G$,
over which the $kG$-modules of finite projective dimension are
precisely those $kG$-modules that have finite projective dimension
over $kH$ for any $\mathfrak{C}$-subgroup $H \subseteq G$; see
\cite[Definition 5.6]{Bi}. Then, the groups of type $\Phi$ over
$k$ are precisely the groups in $\Phi \mathfrak{F}$.

We may consider a variation of the operation $\Phi\_\!\_$ for
flat modules and define for any group class $\mathfrak{C}$ the
class $\Phi_{flat} \mathfrak{C}$ to be the class consisting of
those groups $G$, over which the $kG$-modules of finite flat
dimension are precisely those $kG$-modules having finite flat
dimension over $kH$ for any $\mathfrak{C}$-subgroup $H \subseteq G$.
Starting from the class $\Phi \mathfrak{F}$ of groups of type
$\Phi$ over $k$ (which contains the class $\mathfrak{F}$ of
finite groups), we define the classes $\mathfrak{X}_{\alpha}$
for any ordinal $\alpha$, using transfinite induction, as follows:

(i) $\mathfrak{X}_0 = \Phi \mathfrak{F}$,

(ii) $\mathfrak{X}_{\alpha +1} =
      {\scriptstyle{{\bf LH}}}\mathfrak{X}_{\alpha} \cup
      \Phi_{flat} \mathfrak{X}_{\alpha}$
     for any ordinal $\alpha$ and

(iii) $\mathfrak{X}_{\alpha} =
       \bigcup_{\beta < \alpha}\mathfrak{X}_{\beta}$
      for any limit ordinal $\alpha$.
\newline
Let $\mathfrak{X}$ be the class consisting of those groups which
are contained in $\mathfrak{X}_{\alpha}$, for some ordinal $\alpha$.
Then, $\mathfrak{X}$ contains all
${\scriptstyle{{\bf LH}}}\mathfrak{F}$-groups and all groups of type
$\Phi$ over $k$. In fact, $\mathfrak{X}$ is the smallest class of
groups that contains all groups of type $\Phi$ over $k$ and is
closed under both operations ${\scriptstyle{{\bf LH}}}\_\!\_$ and
$\Phi_{flat}\_\!\_$ (i.e.\
${\scriptstyle{{\bf LH}}}\mathfrak{X} = \mathfrak{X}$ and
$\Phi_{flat} \mathfrak{X} = \mathfrak{X}$.)

Analogously, we may consider a variation of the operation
$\Phi\_\!\_$ for injective modules and define for any group
class $\mathfrak{C}$ the class $\Phi_{inj} \mathfrak{C}$ to
be the class consisting of those groups $G$, over which the
$kG$-modules of finite injective dimension are precisely those
$kG$-modules having finite injective dimension over $kH$ for any
$\mathfrak{C}$-subgroup $H \subseteq G$. Starting again from the
class $\Phi \mathfrak{F}$ of groups of type $\Phi$ over $k$, we
define the classes $\mathfrak{Y}_{\alpha}$ for any ordinal $\alpha$,
using transfinite induction, as follows:

(i) $\mathfrak{Y}_0 = \Phi \mathfrak{F}$,

(ii) $\mathfrak{Y}_{\alpha +1} =
      {\scriptstyle{{\bf LH}}}\mathfrak{Y}_{\alpha} \cup
      \Phi_{inj} \mathfrak{Y}_{\alpha}$
     for any ordinal $\alpha$ and

(iii) $\mathfrak{Y}_{\alpha} =
       \bigcup_{\beta < \alpha}\mathfrak{Y}_{\beta}$
      for any limit ordinal $\alpha$.
\newline
Let $\mathfrak{Y}$ be the class consisting of those groups
which are contained in $\mathfrak{Y}_{\alpha}$, for some
ordinal $\alpha$. Then, $\mathfrak{Y}$ contains all
${\scriptstyle{{\bf LH}}}\mathfrak{F}$-groups and all groups
of type $\Phi$ over $k$. In fact, $\mathfrak{Y}$ is the smallest
class of groups that contains all groups of type $\Phi$ over $k$
and is closed under both operations ${\scriptstyle{{\bf LH}}}\_\!\_$
and $\Phi_{inj}\_\!\_$ (i.e.\
${\scriptstyle{{\bf LH}}}\mathfrak{Y} = \mathfrak{Y}$ and
$\Phi_{inj} \mathfrak{Y} = \mathfrak{Y}$.)

\vspace{0.1in}

\noindent
{\sc III.\ Cotorsion pairs and Hovey triples.}
Let ${\mathcal E}$ be an additive full and extension-closed
subcategory of the category of modules over a ring $R$; then,
$\mathcal{E}$ is an exact category in the sense of \cite{Q}.
The ${\rm Ext}^1$-pairing induces an orthogonality relation
between subclasses of ${\mathcal E}$. If
${\mathcal A} \subseteq {\mathcal E}$, then the left
Ext$^1$-orthogonal $^{\perp}{\mathcal A}$ of ${\mathcal A}$
is the class consisting of those modules $E \in {\mathcal E}$,
which are such that ${\rm Ext}^1_R(E,A)=0$ for all modules
$A \in {\mathcal A}$. Analogously, the right Ext$^1$-orthogonal
${\mathcal A}^{\perp}$ of ${\mathcal A}$ is the class consisting
of those modules $E' \in {\mathcal E}$, which are such that
${\rm Ext}^1_R(A,E')=0$ for all modules $A \in {\mathcal A}$.

Let ${\mathcal C},{\mathcal D}$ be two subclasses of ${\mathcal E}$.
Then, the pair $({\mathcal C},{\mathcal D})$ is a cotorsion pair in
${\mathcal E}$ (cf.\ \cite[Definition 7.1.2]{EJ}) if
${\mathcal C} = {^{\perp} {\mathcal D}}$ and
${\mathcal C} ^{\perp} = {\mathcal D}$. The cotorsion pair is
called hereditary if ${\rm Ext}^i_R(C,D)=0$ for all $i>0$ and
all modules $C \in \mathcal{C}$ and $D \in \mathcal{D}$. If
$\mathcal{C}$ (resp.\ $\mathcal{D}$) contains all projective (resp.\
injective) modules, then the latter condition is equivalent to the
assertion that $\mathcal{C}$ (resp.\ $\mathcal{D}$) is closed under
kernels of epimorphisms (resp.\ under cokernels of monomorphisms).
The cotorsion pair is complete if for any $E \in {\mathcal E}$
there exist two short exact sequences of modules, usually called
approximation sequences
\[ 0 \longrightarrow D \longrightarrow C \longrightarrow E
     \longrightarrow 0
   \;\;\; \mbox{and } \;\;\;
   0 \longrightarrow E \longrightarrow D' \longrightarrow C'
     \longrightarrow 0 , \]
where $C,C' \in {\mathcal C}$ and $D,D' \in {\mathcal D}$. We
say that the cotorsion pair is cogenerated by a set {\tt S} of
modules if $\mathcal{D} = {\tt S}^{\perp}$ (and
$\mathcal{C} = \, \! ^{\perp}({\tt S} ^{\perp})$). Assuming that
$\mathcal{E}$ is the full category of $R$-modules, any cotorsion
pair which is cogenerated by a set is complete; this result is
proved by Eklof and Trlifaj in \cite[Theorem 10]{ETz}.

A Hovey triple on ${\mathcal E}$ is a triple
$({\mathcal C} , {\mathcal W} , {\mathcal F})$ of subclasses of
${\mathcal E}$, which are such that the pairs
$({\mathcal C} , {\mathcal W} \cap {\mathcal F})$ and
$({\mathcal C} \cap {\mathcal W} , {\mathcal F})$ are complete
cotorsion pairs in $\mathcal{E}$ and the class ${\mathcal W}$ is
closed under direct summands and satisfies the 2-out-of-3 property
for short exact sequences in ${\mathcal E}$. Based on the work of
Hovey \cite{Hov2}, Gillespie \cite{G} has established a bijection
between Hovey triples and exact model structures on ${\mathcal E}$,
in the case where $\mathcal{E}$ is closed under direct summands;
cf.\ \cite[Theorem 3.3]{G}. In the context of the latter bijection,
it is proved in \cite[Proposition 5.2]{G} that for an exact model
structure on ${\mathcal E}$ whose associated complete cotorsion
pairs are hereditary, the class ${\mathcal C} \cap {\mathcal F}$
is a Frobenius category with projective-injective objects equal
to ${\mathcal C} \cap {\mathcal W} \cap {\mathcal F}$. A result of
Happel \cite{Ha} then implies that the associated stable category,
i.e.\ ${\mathcal C} \cap {\mathcal F}$ modulo
${\mathcal C} \cap {\mathcal W} \cap {\mathcal F}$, is triangulated.
In view of \cite[Proposition 4.4 and Corollary 4.8]{G}, we conclude
that the homotopy category of an exact model structure as above is
triangulated equivalent to the stable category of the Frobenius
category ${\mathcal C} \cap {\mathcal F}$.

\section{Cofibrant-flat modules}

\noindent
In this section, we consider the class of cofibrant-flat modules
and obtain some basic properties of that class. We relate these
modules to Gorenstein flat modules and show that the class of
cofibrant-flat modules is always the left hand side of a complete
hereditary cotorsion pair.

\vspace{0.1in}

\noindent
{\sc I.\ Cofibrant-flat and Gorenstein flat modules.}
A $kG$-module $M$ is cofibrant-flat if the (diagonal) $kG$-module
$M \otimes_k B$ is flat; we denote by ${\tt Cof.flat}(kG)$ the class
of cofibrant-flat modules. It is clear that any cofibrant module is
cofibrant-flat. Since the $kG$-module $B$ is $k$-free (and hence
$k$-flat), any flat $kG$-module is cofibrant-flat. It follows that
the following conditions are equivalent for a $kG$-module $N$:

(i) There is a non-negative integer $n$ and an exact sequence of
$kG$-modules
\[ 0 \longrightarrow M_n \longrightarrow \cdots \longrightarrow M_1
     \longrightarrow M_0 \longrightarrow N \longrightarrow 0 , \]
such that $M_0,M_1, \ldots ,M_n \in {\tt Cof.flat}(kG)$.

(ii) $\mbox{fd}_{kG}(N \otimes_k B) < \infty$.
\newline
If these conditions are satisfied, we say that $N$ has finite
cofibrant-flat dimension; we denote by $\overline{\tt Cof.flat}(kG)$
the class of these $kG$-modules. It is clear that
${\tt Cof.flat}(kG) \subseteq \overline{\tt Cof.flat}(kG)$. The
following criterion for a $\overline{\tt Cof.flat}(kG)$-module
to be cofibrant-flat is a variant of a result by Benson; cf.\
\cite[Lemma 4.5(ii)]{Ben}.

\begin{Lemma}\label{lem:Ben.flat}
Let $M \in \overline{\tt Cof.flat}(kG)$ and assume that
${\rm Ext}^1_{kG}(M,N)=0$ for any cotorsion $kG$-module $N$
of finite flat dimension. Then, $M \in {\tt Cof.flat}(kG)$.
\end{Lemma}
\vspace{-0.05in}
\noindent
{\em Proof.}
Assume that $M$ is not cofibrant-flat, so that
$n = \mbox{fd}_{kG}(M \otimes_k B)$ is strictly positive. We
consider a special flat precover of $M \otimes_kB$
\[ 0 \longrightarrow K \longrightarrow F
     \stackrel{p}{\longrightarrow} M \otimes_k B
     \longrightarrow 0 , \]
i.e.\ a short exact sequence of $kG$-modules, where $F$ is flat
and $K$ is cotorsion; of course, $K$ has finite flat dimension.
We also consider the commutative diagram with exact rows
\[
\begin{array}{ccccccccc}
 0 & \longrightarrow & L & \longrightarrow
   & F & \longrightarrow & M \otimes_k \overline{B}
   & \longrightarrow & 0 \\
 & & \!\!\! {\scriptstyle{q}} \downarrow
 & & \!\!\! {\scriptstyle{p}} \downarrow
 & & \parallel & & \\
 0 & \longrightarrow & M & \longrightarrow
   & M \otimes_k B & \longrightarrow
   & M \otimes_k \overline{B} & \longrightarrow & 0
\end{array}
\]
where the bottom row is obtained by tensoring with $M$ over
$k$ the $k$-split short exact sequence of $kG$-modules
\[  0 \longrightarrow k \longrightarrow B
      \longrightarrow \overline{B} \longrightarrow 0 . \]
Since $p$ is surjective, an application of the snake lemma
shows that $q$ is also surjective and
$\mbox{ker} \, q = \mbox{ker} \, p = K$. Then, our assumption
on $M$ implies that the epimorphism $q$ splits, so that $M$ is
a direct summand of $L$. It follows that $M \otimes_k B$ is
a direct summand of $L \otimes_k B$. This is a contradiction,
since $M \otimes_k B$ has flat dimension equal to $n$, whereas
the short exact sequence
\[ 0 \longrightarrow L \otimes_k B
     \longrightarrow F \otimes_k B
     \longrightarrow M \otimes_k B \otimes_k \overline{B}
     \longrightarrow 0 \]
shows that $L \otimes _k B$ has flat dimension $\leq n-1$.
\hfill $\Box$

\vspace{0.1in}

\noindent
Cofibrant-flat modules are closely related to Gorenstein flat
modules. Recall that an acyclic complex of flat $kG$-modules
is called totally acyclic if it remains acyclic after applying
the functor $I \otimes_{kG} \_\!\_$ for any injective $kG$-module
$I$. A $kG$-module $M$ is called Gorenstein flat if it is a cokernel
of a totally acyclic complex of flat $kG$-modules. Let ${\tt GFlat}(kG)$
denote the class of Gorenstein flat $kG$-modules. We refer to
\cite{EJ, Hol} for this notion (which may be developed over any
associative ring). Using the homological version of an elegant
construction by Cornick and Kropholler \cite{CK}, it is shown
in \cite[Proposition 5.2(i)]{ET3} that any cofibrant-flat module
is Gorenstein flat\footnote{The proof of \cite[Proposition 5.2(i)]{ET3}
is given in the special case where the coefficient ring is $\mathbb{Z}$.
The same argument though applies verbatim over any coefficient ring $k$.},
so that we always have ${\tt Cof.flat}(kG) \subseteq {\tt GFlat}(kG)$.
In particular, \cite[Corollary 4.12]{SS} implies that the functor
${\rm Ext}^1_{kG}(\_\!\_,N)$ vanishes on cofibrant-flat modules, whenever
$N$ is a cotorsion $kG$-module of finite flat dimension. Hence, the
sufficient condition described in Lemma 2.1 for a $kG$-module to be
cofibrant-flat is also necessary. We now list some basic properties
of the class ${\tt Cof.flat}(kG)$ of cofibrant-flat modules.

\begin{Proposition}\label{prop:cof.flat}
(i)
${\tt Cof.flat}(kG) = \overline{\tt Cof.flat}(kG) \cap
 {\tt GFlat}(kG)$.

(ii) The class ${\tt Cof.flat}(kG)$ is closed under extensions,
kernels of epimorphisms, direct summands and filtered colimits.
It is also closed under pure submodules and pure quotients.

(iii) Any cofibrant-flat module is a cokernel of an acyclic complex
of flat $kG$-modules, all of whose cokernels are cofibrant-flat.

(iv) Any cotorsion cofibrant-flat module is a cokernel of an
acyclic complex of flat cotorsion $kG$-modules, all of whose
cokernels are cotorsion cofibrant-flat.

\end{Proposition}
\vspace{-0.05in}
\noindent
{\em Proof.}
(i) We know that
${\tt Cof.flat}(kG) \subseteq \overline{\tt Cof.flat}(kG)$
and ${\tt Cof.flat}(kG) \subseteq {\tt GFlat}(kG)$. On the
other hand, if
$M \in \overline{\tt Cof.flat}(kG) \cap {\tt GFlat}(kG)$,
then ${\rm Ext}^1_{kG}(M,N)=0$ for any cotorsion $kG$-module
$N$ of finite flat dimension; cf.\ \cite[Corollary 4.12]{SS}.
Hence, Lemma 2.1 implies that $M$ is cofibrant-flat, as needed.

(ii) Since the class of flat $kG$-modules is closed under
extensions, kernels of epimorphisms, direct summands and
filtered colimits, it follows readily that the class of
cofibrant-flat modules enjoys these closure properties
as well. We now let
\[ 0 \longrightarrow M' \longrightarrow M
     \longrightarrow M'' \longrightarrow 0 \]
be a pure exact sequence of $kG$-modules. It is well-known
that such an exact sequence is the filtered colimit of split
short exact sequences. Since the endofunctor $\_\!\_ \otimes_kB$
of the category of $kG$-modules preserves filtered colimits and
split short exact sequences, the induced short exact sequence of
$kG$-modules
\[ 0 \longrightarrow M' \otimes_kB \longrightarrow M \otimes_kB
     \longrightarrow M'' \otimes_kB \longrightarrow 0 \]
is also a filtered colimit of split short exact sequences; hence,
it is pure exact. We now assume that the $kG$-module $M$ is
cofibrant-flat, so that the $kG$-module $M \otimes_kB$ is flat.
The class of flat $kG$-modules is closed under pure submodules
and pure quotients and hence the $kG$-modules $M' \otimes_kB$
and $M'' \otimes_kB$ are also flat, so that $M',M''$ are both
cofibrant-flat.

(iii) Let $N$ be a cofibrant-flat module. Since all flat
$kG$-modules are cofibrant-flat, it follows from (ii) above that
the cokernels of any flat resolution of $N$ are cofibrant-flat.
This provides us with the left half of the required acyclic
complex of flat $kG$-modules. On the other hand, being
cofibrant-flat, the $kG$-module $N$ is Gorenstein flat. Hence,
there exists an exact sequence of $kG$-modules
\[ 0 \longrightarrow N \longrightarrow F_0 \longrightarrow F_1
     \longrightarrow \cdots , \]
where the $F_i$'s are flat (and hence cofibrant-flat) and all
cokernels are Gorenstein flat. Since these cokernels have also
finite cofibrant-flat dimension, it follows from (i) above that
they are actually cofibrant-flat. We therefore obtain the right
half of the required acyclic complex of flat $kG$-modules.

(iv) The argument is similar to that used in (iii) above.
Let $N$ be a cotorsion cofibrant-flat module. The left half
of the required acyclic complex of flat cotorsion $kG$-modules
is obtained by considering successive flat covers of cotorsion
cofibrant-flat modules, using the closure of ${\tt Cof.flat}(kG)$
under kernels of epimorphisms; cf.\ (ii) above. The right half
of the complex is obtained by using successively
\cite[Lemma 4.1]{E2} and (i) above. \hfill $\Box$

\vspace{0.1in}

\noindent
Let us denote by $\mathcal{F}(kG)$ the class of those
$kG$-modules that appear as cokernels of acyclic complexes
of flat $kG$-modules. It is clear that
${\tt GFlat}(kG) \subseteq \mathcal{F}(kG)$ and equality holds
if and only if any acyclic complex of flat $kG$-modules is
totally acyclic. In complete analogy with the conjecture by
Dembegioti and Talelli \cite{DT} regarding cofibrant modules,
we may ask whether $\mathcal{F}(kG) \subseteq {\tt Cof.flat}(kG)$,
so that ${\tt Cof.flat}(kG) = {\tt GFlat}(kG) = \mathcal{F}(kG)$.
We recall from $\S $1.II that $\mathfrak{X}$ is the smallest
class of groups that contains all groups of type $\Phi$ over
$k$ and satisfies the equalities
${\scriptstyle{{\bf LH}}}\mathfrak{X} = \mathfrak{X}$ and
$\Phi_{flat} \mathfrak{X} = \mathfrak{X}$. The class
$\mathfrak{X}$ contains all
${\scriptstyle{{\bf LH}}}\mathfrak{F}$-groups.

\begin{Proposition}\label{prop:C=GF1}
Assume that any acyclic complex of flat $k$-modules is
pure-acyclic\footnote{This is the case if, for example,
$k$ has finite weak global dimension.} and let $G$ be an $\mathfrak{X}$-group. Then,
${\tt Cof.flat}(kG) = {\tt GFlat}(kG) = \mathcal{F}(kG)$.
\end{Proposition}
\vspace{-0.05in}
\noindent
{\em Proof.}
We shall prove that $\mathcal{F}(kG) \subseteq {\tt Cof.flat}(kG)$.
Denote by $\mathcal{B}(kG)$ the class consisting of those $kG$-modules
$B'$ which are projective over any finite subgroup of $G$ and let
${\tt Cof.flat}'(kG)$ be the class consisting of those $kG$-modules
$M$, for which the diagonal $kG$-module $M \otimes_kB'$ is flat for
all $B' \in \mathcal{B}(kG)$. Since $B \in \mathcal{B}(kG)$, it
follows that ${\tt Cof.flat}'(kG) \subseteq {\tt Cof.flat}(kG)$.
Hence, it suffices to show that
$\mathcal{F}(kG) \subseteq {\tt Cof.flat}'(kG)$. To that end, we
consider the class
\[ \mathfrak{Z} =
   \{ H : \mathcal{F}(kH) \subseteq {\tt Cof.flat}'(kH) \} \]
and show that $\mathfrak{X} \subseteq \mathfrak{Z}$. In view of the
definition of $\mathfrak{X}$, it suffices to show that:
\newline
(i) all groups of type $\Phi$ over $k$ are contained in
$\mathfrak{Z}$ and
\newline
(ii) ${\scriptstyle{{\bf LH}}}\mathfrak{Z} = \mathfrak{Z}$ and
$\Phi_{flat} \mathfrak{Z} = \mathfrak{Z}$.

In order to prove (i), consider a group $H$ of type $\Phi$ over
$k$ and let $M \in \mathcal{F}(kH)$. Since flat $kH$-modules are
also $k$-flat, our assumption on $k$ implies that $M$ is $k$-flat.
We also consider a $kH$-module $B' \in \mathcal{B}(kH)$ and note
that $B'$ has finite projective dimension; in particular, $B'$ has
finite flat dimension, say $\mbox{fd}_{kH}B'  = n$. It follows
readily that $\mbox{fd}_{kH}(M \otimes_kB') \leq \mbox{fd}_{kH}B'=n$.
This is true for any module $M \in \mathcal{F}(kH)$. Since $M$ is an
$n$-th flat syzygy of another $kH$-module $M' \in \mathcal{F}(kH)$
and $B'$ is $k$-flat (being $k$-projective), we conclude that the
$kH$-module $M \otimes _kB'$ is actually flat. It follows that
$M \in {\tt Cof.flat}'(kH)$ and hence $H \in \mathfrak{Z}$.

We shall now prove (ii) and begin by verifying that $\mathfrak{Z}$
is ${\scriptstyle{{\bf LH}}}$-closed. Consider an
${\scriptstyle{{\bf LH}}}\mathfrak{Z}$-group $H$ and let
$M \in \mathcal{F}(kH)$ and $B' \in \mathcal{B}(kH)$. Then, we can
show that the $kH$-module $M \otimes_k B'$ is flat, following verbatim
the arguments in the proof of \cite[Proposition 5.2(ii)]{ET3} and
proving that $M \otimes_k B'$ is flat over any
${\scriptstyle{{\bf LH}}}\mathfrak{Z}$-subgroup of $H$. The first step
of the transfinite induction therein (showing that $M \otimes_k B'$ is
flat over any ${\scriptstyle{{\bf H}}}\mathfrak{Z}$-subgroup of $H$)
follows from the definition of $\mathfrak{Z}$. We conclude that
$M \in {\tt Cof.flat}'(kH)$ and hence $H \in \mathfrak{Z}$.

Finally, we shall prove that $\mathfrak{Z}$ is $\Phi_{flat}$-closed.
To that end, we consider a $\Phi_{flat}\mathfrak{Z}$-group $H$ and
let $M \in \mathcal{F}(kH)$ and $B' \in \mathcal{B}(kH)$. Then, the
definition of $\mathfrak{Z}$ implies that $M \otimes_kB'$ is flat
over any $\mathfrak{Z}$-subgroup of $H$. Since $H$ is a
$\Phi_{flat}\mathfrak{Z}$-group, we conclude that the $kH$-module
$M \otimes_kB'$ has finite flat dimension. This is true for any
$kH$-module $M \in \mathcal{F}(kH)$. In other words, having fixed
$B' \in \mathcal{B}(kH)$, we have established the finiteness of the
flat dimension of $M' \otimes_kB'$ for any $M' \in \mathcal{F}(kH)$.
The class $\mathcal{F}(kH)$ being closed under direct sums, it follows
easily that there is an upper bound, say $n$, for the flat dimensions
$\mbox{fd}_{kH}(M' \otimes_kB')$, $M' \in \mathcal{F}(kH)$. Since $M$
is an $n$-th flat syzygy of another $kG$-module $M' \in \mathcal{F}(kG)$
and $B'$ is $k$-flat, we conclude that the $kG$-module $M \otimes _kB'$
is actually flat. It follows that $M \in {\tt Cof.flat}'(kH)$ and hence
$H \in \mathfrak{Z}$. \hfill $\Box$

\vspace{0.1in}

\noindent
It is not known whether the class of Gorenstein flat modules (over
a general ring) is necessarily closed under pure submodules and pure
quotients, even though \cite[Theorems 2.6 and 3.6]{Hol} imply that
this is the case over right coherent rings; see also \cite[Proposition
4.13]{SS}. We provide another class of rings (which are often
non-coherent) for which this assertion is also true.

\begin{Corollary}
Assume that any acyclic complex of flat $k$-modules is
pure-acyclic\ and $G$ is an $\mathfrak{X}$-group. Then,
the class ${\tt GFlat}(kG)$ is closed under pure submodules
and pure quotients.
\end{Corollary}
\vspace{-0.05in}
\noindent
{\em Proof.}
This is an immediate consequence of Proposition 2.2(ii) and Proposition
2.3. \hfill $\Box$

\vspace{0.1in}

\noindent
{\bf Remark 2.5.}
Assuming that ${\tt Cof.flat}(kG) = {\tt GFlat}(kG)$, Lemma
\ref{lem:Ben.flat} states that, within the class of $kG$-modules
of finite Gorenstein flat dimension, the Gorenstein flat $kG$-modules
are precisely those $kG$-modules which are left Ext$^1$-orthogonal
to the cotorsion $kG$-modules of finite flat dimension; see
\cite[Lemma 5.4]{C6}. In the same way, the equality
${\tt Cof.flat}(kG) = {\tt GFlat}(kG)$ reduces the assertions stated
in the first sentence of Proposition \ref{prop:cof.flat}(ii) to certain
well-known properties of Gorenstein flat modules (cf.\ \cite[Corollary
4.12]{SS}), whereas Proposition \ref{prop:cof.flat}(iii) is essentially
a consequence of the very definition of Gorenstein flat modules.
\addtocounter{Lemma}{1}

\vspace{0.1in}

\noindent
Let $M$ be a cofibrant module and invoke
\cite[Proposition 2.2(ii)]{ER} to express $M$ as a cokernel
of an acyclic complex of projective $kG$-modules, say $P$,
all of whose cokernels are cofibrant. Since
${\tt Cof}(kG) \subseteq {\tt Cof.flat}(kG) \subseteq {\tt GFlat}(kG)$,
these cokernels are Gorenstein flat and hence the complex
$P \otimes_{kG}I$ is acyclic for any injective $kG$-module $I$.
The cokernels of those acyclic complexes of projective $kG$-modules
that remain acyclic by applying the functor $\_\!\_ \otimes_{kG}I$
for any injective $kG$-module $I$ are called projectively coresolved
Gorenstein flat; cf.\ \cite[$\S 4$]{SS}. We denote by ${\tt PGF}(kG)$
the class of these $kG$-modules. It follows that $M$ is projectively
coresolved Gorenstein flat, and hence
${\tt Cof}(kG) \subseteq {\tt PGF}(kG)$; cf.\ \cite[Proposition 8.2]{S}.

\begin{Corollary}
(i) ${\tt Cof.flat}(kG) \cap {\tt PGF}(kG) = {\tt Cof}(kG)$.

(ii) If all Gorenstein projective $kG$-modules are Gorenstein
flat, then
${\tt Cof.flat}(kG) \cap {\tt GProj}(kG) = {\tt Cof}(kG)$.
\end{Corollary}
\vspace{-0.05in}
\noindent
{\em Proof.}
(i) We know that
${\tt Cof}(kG) \subseteq {\tt Cof.flat}(kG) \cap {\tt PGF}(kG)$.
To show the reverse inclusion, we consider a cofibrant-flat
projectively coresolved Gorenstein flat $kG$-module $M$ and let
$P$ be an acyclic complex of projective $kG$-modules with $M=C_0P$,
whose cokernels $(C_nP)_n$ are projectively coresolved Gorenstein
flat. As these cokernels have finite cofibrant-flat dimension,
Proposition 2.2(i) implies that they are all cofibrant-flat. It
follows that $P \otimes_kB$ is an acyclic complex of projective
$kG$-modules, whose cokernels $(C_nP \otimes_kB)_n$ are flat.
Invoking \cite[Remark 2.15]{N2}, we conclude that these cokernels
are projective. In particular, $M \otimes_kB = C_0P \otimes_kB$
is projective, so that $M$ is cofibrant.

(ii) If all Gorenstein projective $kG$-modules are Gorenstein
flat, then all Gorenstein projective $kG$-modules are projectively
coresolved Gorenstein flat, so that ${\tt GProj}(kG) = {\tt PGF}(kG)$.
Hence, the result follows from (i) above. \hfill $\Box$

\vspace{0.1in}

\noindent
{\sc II.\ The cotorsion pair
$({\tt Cof.flat}(kG),{\tt Cof.flat}(kG)^{\perp})$.}
Let $\lambda$ be a cardinal number. We recall from \cite[Definition 2.1]{ELR}
that a class $\mathcal{K}$ of $kG$-modules is called $\lambda$-Kaplansky
if for any module $M \in \mathcal{K}$ and any element $x \in M$ there
exists a submodule $N \subseteq M$, such that $x \in N$,
$\mbox{card} \, N \leq \lambda$ and $N,M/N \in \mathcal{K}$. We say
that $\mathcal{K}$ is a Kaplansky class if it is $\lambda$-Kaplansky,
for some cardinal number $\lambda$. As shown by Holm and J\o rgensen
in \cite[Proposition 3.2]{HJ}, any class of $kG$-modules, which is
closed under pure submodules and pure quotients, is Kaplansky. In
particular, Lemma 2.2(ii) implies that the class of cofibrant-flat
modules is Kaplansky. Using this property of ${\tt Cof.flat}(kG)$,
we may construct filtrations of cofibrant-flat modules, whose
successive quotients are cofibrant-flat modules of controlled size.

\begin{Proposition}\label{prop:caf}
There is a cardinal number $\lambda$, such that any cofibrant-flat
module admits a continuous ascending filtration by $\lambda$-generated
cofibrant-flat modules.
\end{Proposition}
\vspace{-0.05in}
\noindent
{\em Proof.}
Let $\lambda$ be a cardinal number, such that the class
${\tt Cof.flat}(kG)$ is $\lambda$-Kaplansky. We consider a cofibrant-flat
module $M$ and choose a well-ordered set of generators
$(x_{\alpha})_{\alpha < \tau}$ of it. We shall construct by transfinite
induction a continuous ascending family of $kG$-submodules
$(M_{\alpha})_{\alpha \leq \tau}$ of $M$, such that:

(i) $M_0=0$,

(ii) $x_{\alpha} \in M_{\alpha +1}$ for any $\alpha < \tau$,

(iii) $M/M_{\alpha}$ is cofibrant-flat for all $\alpha < \tau$ and

(iv) $M_{\alpha +1}/M_{\alpha}$ is cofibrant-flat and
$\lambda$-generated for any $\alpha < \tau$.
\newline
We begin the induction, by letting $M_0=0$. For the inductive
step, we assume that $\beta \leq \tau$ is an ordinal and the
$M_{\alpha}$'s have been constructed for all $\alpha < \beta$,
so that properties (i)-(iv) above hold. In the case where
$\beta$ is a limit ordinal, we define
$M_{\beta} = \bigcup_{\alpha < \beta} M_{\alpha}$. In this way,
property (iii) still holds, since $M/M_{\beta}$ is the filtered
colimit of the family $(M/M_{\alpha})_{\alpha < \beta}$ of
cofibrant-flat modules and ${\tt Cof.flat}(kG)$ is closed under
filtered colimits, in view of Proposition 2.2(ii). We now consider
the case where $\beta = \alpha +1$ is a successor ordinal. Since
$M/M_{\alpha}$ is cofibrant-flat and the class ${\tt Cof.flat}(kG)$
is $\lambda$-Kaplansky, there exists a $kG$-submodule
$M_{\alpha +1} \subseteq M$ containing $M_{\alpha}$ and
$x_{\alpha}$, such $M_{\alpha +1}/M_{\alpha}$ is a
$\lambda$-generated cofibrant-flat module and the quotient
$(M/M_{\alpha})/(M_{\alpha +1}/M_{\alpha}) = M/M_{\alpha +1}$
is cofibrant-flat. The inductive step of the construction is
thus completed. Since the $kG$-submodule $M_{\tau}$ contains all
of the generators $(x_{\alpha})_{\alpha < \tau}$ of $M$, it follows
that $M_{\tau}=M$. Hence, $(M_{\alpha})_{\alpha \leq \tau}$ is the
required continuous ascending filtration of $M$ by
$\lambda$-generated cofibrant-flat modules. \hfill $\Box$

\begin{Theorem}\label{thm:cp-flat}
The pair
$\left( {\tt Cof.flat}(kG),{\tt Cof.flat}(kG)^{\perp} \right)$
is a hereditary cotorsion pair in the category of $kG$-modules,
which is cogenerated by a set. Its kernel is the class of
flat cotorsion $kG$-modules.
\end{Theorem}
\vspace{-0.05in}
\noindent
{\em Proof.}
In view of Proposition 2.7, there exists a cardinal number
$\lambda$, such that any cofibrant-flat module admits a continuous
ascending filtration by $\lambda$-generated cofibrant-flat modules.
We let ${\mathcal S}$ be a {\em set} of representatives of the
isomorphism classes of all $\lambda$-generated cofibrant-flat
modules and prove that
\[ \left( {\tt Cof.flat}(kG),{\tt Cof.flat}(kG)^{\perp} \right) =
   \left( ^{\perp} \! \left( \mathcal{S}^{\perp} \right) \! ,
   \mathcal{S}^{\perp} \right) \]
is the cotorsion pair cogenerated by $\mathcal{S}$. Of course,
it only suffices to show that
${\tt Cof.flat}(kG) =
 \, \! ^{\perp} \! \left( \mathcal{S}^{\perp} \right)$.
In view of \cite[Corollary 7.3.5]{EJ}, the class
$^{\perp} \! \left( \mathcal{S}^{\perp} \right)$ consists of the
direct summands of those $kG$-modules that admit a continuous
ascending filtration by modules in $\mathcal{S}$. Hence, the
inclusion
${\tt Cof.flat}(kG) \subseteq \, \!
 ^{\perp} \! \left( \mathcal{S}^{\perp} \right)$
is an immediate consequence of Proposition 2.7. Since $\mathcal{S}$
is contained in ${\tt Cof.flat}(kG)$ and the latter class is closed
under direct summands, the reverse inclusion will follow if we show
that any $kG$-module that admits a continuous ascending filtration
by ${\tt Cof.flat}(kG)$-modules is also contained in
${\tt Cof.flat}(kG)$. The latter claim follows readily by transfinite
induction on the length of the filtration, since the class of
cofibrant-flat modules is closed under extensions and filtered unions;
cf.\ Proposition 2.2(ii).

The cotorsion pair
$\left( {\tt Cof.flat}(kG),{\tt Cof.flat}(kG)^{\perp} \right)$ is
hereditary, since the left hand class ${\tt Cof.flat}(kG)$ is
closed under kernels of epimorphisms; cf.\ Proposition 2.2(ii). It
only remains to show that the kernel of the cotorsion pair is the
class of flat cotorsion modules, i.e.\ that
\[ {\tt Cof.flat}(kG) \cap {\tt Cof.flat}(kG)^{\perp} =
   {\tt Flat}(kG) \cap {\tt Cotor}(kG) . \]
Let $M$ be a $kG$-module contained in the intersection
${\tt Cof.flat}(kG) \cap {\tt Cof.flat}(kG)^{\perp}$. Then, the
functor ${\rm Ext}_{kG}^1(\_\!\_,M)$ vanishes on all cofibrant-flat
modules and, in particular, on all flat $kG$-modules; hence, $M$ is
cotorsion. Since $M$ is cofibrant-flat, the $kG$-module
$M \otimes_k\overline{B}$ is also cofibrant-flat. We conclude that
the short exact sequence of $kG$-modules
\[  0 \longrightarrow M \longrightarrow M \otimes_kB
      \longrightarrow M \otimes_k\overline{B} \longrightarrow 0 , \]
which is obtained by tensoring with $M$ over $k$ the $k$-split short
exact sequence of $kG$-modules
\[  0 \longrightarrow k \longrightarrow B
      \longrightarrow \overline{B} \longrightarrow 0 , \]
splits. Therefore, $M$ is a direct summand of the flat $kG$-module
$M \otimes_kB$ and hence it is itself flat. For the reverse inclusion,
we note that any flat cotorsion $kG$-module is flat and hence
cofibrant-flat, whereas
\[ {\tt Flat}(kG) \cap {\tt Cotor}(kG)  \subseteq
   {\tt GFlat}(kG)^{\perp} \subseteq {\tt Cof.flat}(kG)^{\perp} . \]
Here, the first inclusion follows from \cite[Corollary 4.12]{SS} and
the second one is obvious, since any cofibrant-flat module is Gorenstein
flat. \hfill $\Box$

\vspace{0.1in}

\noindent
We now list a few properties of the class ${\tt Cof.flat}(kG)^{\perp}$.
Since any flat $kG$-module is cofibrant-flat, all $kG$-modules in
${\tt Cof.flat}(kG)^{\perp}$ are cotorsion.

\begin{Proposition}\label{prop:cof+.flat}
The class ${\tt Cof.flat}(kG)^{\perp}$ is closed under direct summands
and has the 2-out-of-3 property for short exact sequences of cotorsion
$kG$-modules. It contains all cotorsion $kG$-modules of finite flat or
injective dimension and all $kG$-modules of the form ${\rm Hom}_k(B,N)$,
where $N$ is a cotorsion $kG$-module.
\end{Proposition}
\vspace{-0.05in}
\noindent
{\em Proof.}
The class ${\tt Cof.flat}(kG)^{\perp}$ is obviously closed under
direct summands and extensions. Since the cotorsion pair
$\left( {\tt Cof.flat}(kG),{\tt Cof.flat}(kG)^{\perp} \right)$
is hereditary, ${\tt Cof.flat}(kG)^{\perp}$ is also closed under
cokernels of monomorphisms. Consider a short exact sequence of
$kG$-modules
\begin{equation}
 0 \longrightarrow N' \longrightarrow N \longrightarrow N''
   \longrightarrow 0 ,
\end{equation}
where $N,N'' \in {\tt Cof.flat}(kG)^{\perp}$ and $N'$ is cotorsion.
We shall prove that $N' \in {\tt Cof.flat}(kG)^{\perp}$. To that
end, let $M$ be a cofibrant-flat module and invoke Proposition
\ref{prop:cof.flat}(iii), in order to find a short exact sequence
of $kG$-modules
\begin{equation}
 0 \longrightarrow M \longrightarrow F \longrightarrow M'
   \longrightarrow 0 ,
\end{equation}
where $F$ is flat and $M'$ is cofibrant-flat. Since both groups
${\rm Ext}^1_{kG}(M',N'')$ and ${\rm Ext}^2_{kG}(M',N)$ are
trivial, it follows from (1) that ${\rm Ext}^2_{kG}(M',N')=0$.
Then, using our assumption that $N'$ is cotorsion, we may conclude
from (2) that ${\rm Ext}^1_{kG}(M,N') = {\rm Ext}^2_{kG}(M',N')$
is the trivial group. Hence, $N' \in {\tt Cof.flat}(kG)^{\perp}$,
as needed.

As we noted in the discussion preceding Proposition 2.2, the class
${\tt Cof.flat}(kG)^{\perp}$ contains all cotorsion $kG$-modules of
finite flat dimension. On the other hand, ${\tt Cof.flat}(kG)^{\perp}$
contains all injective $kG$-modules as well. Hence, an inductive argument
and the 2-out-of-3 property for short exact sequences of cotorsion
$kG$-modules established above show that ${\tt Cof.flat}(kG)^{\perp}$
contains all cotorsion $kG$-modules of finite injective dimension.
Finally, if $N$ is any cotorsion $kG$-module, the Hom-tensor adjunction
for diagonal $kG$-modules implies that the functor
\[ {\rm Ext}^1_{kG}(\_\!\_,\mbox{Hom}_k(B,N)) =
   {\rm Ext}^1_{kG}(\_\!\_ \otimes_k B,N) \]
vanishes on all cofibrant-flat modules and hence
$\mbox{Hom}_k(B,N) \in {\tt Cof.flat}(kG)^{\perp}$. \hfill $\Box$

\section{Approximations of cofibrant-flat modules}

\noindent
In this section, we show that any cofibrant-flat module can be
approximated by cofibrant or flat modules, in a way analogous
to the approximations of Gorenstein flat modules obtained by
\v{S}aroch and \v{S}t$\!$'$\!$ov\'{i}\v{c}ek in \cite{SS}. We
also show that these approximations may be interpreted in terms
of the existence of certain model structures.

\vspace{0.1in}

\noindent
{\sc I.\ Approximations by cofibrant or flat modules.}
Neeman proved in \cite{N2} that, over any ring, the embedding
of the homotopy category of projective modules into the homotopy
category of flat modules admits a right adjoint, whose kernel
consists precisely of the pure-acyclic complexes of flat modules.
Hence, any complex of flat modules $F$ admits a chain map $f$
from a complex of projective modules $Q$, such that the cone
of $f$ is pure-acyclic. It follows that there exists a short
exact sequence of chain complexes
\begin{equation}
 0 \longrightarrow F \longrightarrow Z \longrightarrow P
   \longrightarrow 0 ,
\end{equation}
where $P=Q[1]$ is a complex of projective modules and
$Z = \mbox{cone}(f)$ is a pure-acyclic complex of flat modules.
Of course, if $F$ is acyclic, then $P$ is also acyclic. We shall
use the existence of the short exact sequence (3), in the special
case where the ground ring is the group algebra $kG$, to obtain
a version of \cite[Theorem 4.11]{SS} for $kG$-modules.

\begin{Theorem}
The following conditions are equivalent for a $kG$-module $M$:

(i) $M$ is cofibrant-flat.

(ii) There exists a short exact sequence of $kG$-modules
\[ 0 \longrightarrow M \longrightarrow F \longrightarrow N
     \longrightarrow 0 , \]
where $F$ is flat and $N$ is cofibrant.

(iii) There exists a short exact sequence of $kG$-modules
\[ 0 \longrightarrow K \longrightarrow L \longrightarrow M
     \longrightarrow 0 , \]
where $K$ is flat, $L$ is cofibrant and all $kG$-linear maps
from $K$ to cotorsion $kG$-modules extend to $L$.

(iv) ${\rm Ext}^1_{kG}(M,C)=0$ for any cotorsion $kG$-module
$C \in {\tt Cof}(kG)^{\perp}$.
\end{Theorem}
\vspace{-0.05in}
\noindent
{\em Proof.}
(i)$\rightarrow$(ii): Assume that $M$ is cofibrant-flat and
fix an acyclic complex of flat $kG$-modules $F$, such that
$M = C_0F$ and all cokernels $(C_nF)_n$ are cofibrant-flat;
cf.\ Proposition 2.2(iii). We also fix a short exact sequence
of chain complexes of $kG$-modules as in (3), where $P$ is an
acyclic complex of projective $kG$-modules and $Z$ is a
pure-acyclic complex of flat $kG$-modules. For any integer $n$
we consider the induced short exact sequence of $kG$-modules
\[  0 \longrightarrow C_nF \longrightarrow C_nZ
      \longrightarrow C_nP \longrightarrow 0 . \]
We note that $C_nF$ is cofibrant-flat and $C_nZ$ is flat for
all $n$. Tensoring the above short exact sequence with $B$
over $k$, we obtain a short exact sequence of $kG$-modules
\[  0 \longrightarrow C_nF \otimes_kB
      \longrightarrow C_nZ \otimes_kB
      \longrightarrow C_nP \otimes_kB
      \longrightarrow 0 , \]
so that $\mbox{fd}_{kG}(C_nP \otimes_kB) \leq 1$ for all $n$.
Hence, $P \otimes_kB$ is an acyclic complex of projective
$kG$-modules, all of whose cokernels $(C_nP \otimes_kB)_n$ have
flat dimension $\leq 1$. It follows readily that these cokernels
are flat; in fact, invoking \cite[Remark 2.15]{N2}, we conclude
that these cokernels are actually projective and hence $C_nP$ is
cofibrant for all $n$. The short exact sequence
\[ 0 \longrightarrow C_0F \longrightarrow C_0Z \longrightarrow C_0P
     \longrightarrow 0 \]
is therefore of the type required in (ii).

(ii)$\rightarrow$(iii): We consider a short exact sequence as in
(ii) and fix a short exact sequence of $kG$-modules
\[ 0 \longrightarrow K \longrightarrow P \longrightarrow F
     \longrightarrow 0 , \]
where $P$ is projective and $K$ is, of course, flat. We consider
the pullback of $P \longrightarrow F$ along the monomorphism
$M \longrightarrow F$ and obtain the commutative diagram with
exact rows
\begin{equation}
 \begin{array}{ccccccccc}
  & & 0 & & 0 & & & & \\
  & & \downarrow & & \downarrow & & & & \\
  & & K & = & K & & & & \\
  & & \downarrow & & \downarrow & & & & \\
  0 & \longrightarrow & L & \longrightarrow & P
    & \longrightarrow & N & \longrightarrow & 0 \\
  & & \downarrow & & \downarrow & & \parallel & & \\
  0 & \longrightarrow & M & \longrightarrow & F
    & \longrightarrow & N & \longrightarrow & 0 \\
  & & \downarrow & & \downarrow & & & & \\
  & & 0 & & 0 & & & &
 \end{array}
\end{equation}
Since both $P$ and $N$ are cofibrant, the closure of ${\tt Cof}(kG)$
under kernels of epimorphisms (cf.\ \cite[Proposition 2.2(ii)]{ER})
implies that $L$ is also cofibrant. We claim that the leftmost
vertical short exact sequence is of the type required in (iii).
Indeed, if $C$ is a cotorsion $kG$-module, then any $kG$-linear
map $K \longrightarrow C$ extends to $P$, since $F$ is flat;
in particular, it extends to $L$.

(iii)$\rightarrow$(iv): We consider a short exact sequence
as in (iii) and fix a cotorsion $kG$-module $C$. In view of
our assumption, the additive map
${\rm Ext}^1_{kG}(M,C) \longrightarrow {\rm Ext}^1_{kG}(L,C)$
is injective. If, in addition, $C \in {\tt Cof}(kG)^{\perp}$,
then the group ${\rm Ext}^1_{kG}(L,C)$ is trivial and hence
${\rm Ext}^1_{kG}(M,C)=0$.

(iv)$\rightarrow$(i): Since the class of cofibrant-flat modules
contains all flat and all cofibrant modules, it is clear that
${\tt Cof.flat}(kG)^{\perp} \subseteq
 {\tt Cotor}(kG) \cap {\tt Cof}(kG)^{\perp}$.
Hence, assertion (iv) implies that ${\rm Ext}^1_{kG}(M,C)=0$
for any $C \in {\tt Cof.flat}(kG)^{\perp}$, so that the $kG$-module
$M$ is cofibrant-flat, in view of Theorem 2.8. \hfill $\Box$

\begin{Corollary}
There is an equality
${\tt Cof.flat}(kG)^{\perp} =
 {\tt Cotor}(kG) \cap {\tt Cof}(kG)^{\perp}$,
so that the cotorsion pair
$\left( {\tt Cof.flat}(kG),{\tt Cof.flat}(kG)^{\perp} \right)$
of Theorem 2.8 is the supremum of the cotorsion pairs
$\left( {\tt Flat}(kG),{\tt Cotor}(kG) \right)$ and
$\left( {\tt Cof}(kG),{\tt Cof}(kG)^{\perp} \right)$; cf.\
\cite[Theorem 3.3]{ER}.
\end{Corollary}
\vspace{-0.05in}
\noindent
{\em Proof.}
As we noted during the proof of the implication (iv)$\rightarrow$(i)
in Theorem 3.1, the inclusion
${\tt Cof.flat}(kG)^{\perp} \subseteq
 {\tt Cotor}(kG) \cap {\tt Cof}(kG)^{\perp}$
is obvious. The reverse inclusion follows from the implication
(i)$\rightarrow$(iv) in Theorem 3.1. \hfill $\Box$

\vspace{0.1in}

\noindent
Since cofibrant modules are projectively coresolved Gorenstein
flat, \cite[Theorem 4.4]{SS} implies that ${\rm Ext}^1_{kG}(M,N)=0$
if $M$ is a cofibrant and $N$ is a flat $kG$-module. In other words,
there are inclusions ${\tt Flat}(kG) \subseteq {\tt Cof}(kG)^{\perp}$
and ${\tt Cof}(kG) \subseteq \, \! ^{\perp}{\tt Flat}(kG)$.

\begin{Corollary}
Let $M$ be a cofibrant-flat module.

(i) $M$ is cofibrant if and only if ${\rm Ext}^1_{kG}(M,K)=0$
for any flat $kG$-module $K$.

(ii) $M$ is flat if and only if ${\rm Ext}^1_{kG}(N,M)=0$ for
any cofibrant module $N$.
\end{Corollary}
\vspace{-0.05in}
\noindent
{\em Proof.}
(i) If $M$ is cofibrant, then the group ${\rm Ext}^1_{kG}(M,K)$
is trivial for any flat $kG$-module $K$, in view of the discussion
above. The converse follows from Theorem 3.1(iii).

(ii) If $M$ is flat, then the group ${\rm Ext}^1_{kG}(N,M)$ is
trivial for any cofibrant module $N$, in view of the discussion
above. The converse follows from Theorem 3.1(ii). \hfill $\Box$

\vspace{0.1in}

\noindent
{\bf Remarks 3.4.}
(i) If $M$ is a cofibrant-flat module, the commutative diagram (4)
implies that there exists a short exact sequence of $kG$-modules
\[ 0 \longrightarrow L \longrightarrow M \oplus P
     \longrightarrow F \longrightarrow 0 , \]
where $L$ is cofibrant, $P$ is projective and $F$ is flat. It follows
that ${\tt Cof.flat}(kG)$ is the smallest class of $kG$-modules,
which is closed under extensions and direct summands and
contains all cofibrant and all flat $kG$-modules.

(ii) Regarding the schematic diagram below
\[
\begin{array}{ccccc}
 & & \!\!\! {\tt Cof.flat}(kG) \!\!\! & & \\
 & \!\!\! \nearrow \!\!\! & & \!\!\! \nwarrow \!\!\! & \\
 {\tt Cof}(kG) \!\!\! & & & & \!\!\! {\tt Flat}(kG)  \\
 & \!\!\! \nwarrow \!\!\! & & \!\!\! \nearrow \!\!\! & \\
 & & \!\!\! {\tt Proj}(kG) \!\!\! & &
\end{array}
\]
where all arrows are inclusions, we note that Corollary 3.3
implies that
\[ {\tt Cof.flat}(kG) \cap {\tt Cof}(kG)^{\perp} =
   {\tt Flat}(kG) \;\;\; \mbox{and} \;\;\;
   {\tt Cof.flat}(kG) \cap \, \! ^{\perp}{\tt Flat}(kG) =
   {\tt Cof}(kG) . \]
We also note that the intersection
${\tt Cof}(kG) \cap {\tt Flat}(kG)$ is equal to the class
${\tt Proj}(kG)$ of projective $kG$-modules. Indeed, any
flat and projectively coresolved Gorenstein flat $kG$-module
is necessarily projective; cf.\ \cite[Theorem 4.4]{SS}.
\addtocounter{Lemma}{1}

\begin{Corollary}
The class of cofibrant-flat modules is precisely the class
of pure quotients of cofibrant modules.
\end{Corollary}
\vspace{-0.05in}
\noindent
{\em Proof.}
Any pure quotient of a cofibrant module (in fact, any pure
quotient of a cofibrant-flat module) is cofibrant-flat; cf.\
Proposition 2.2(ii). Conversely, assume that $M$ is cofibrant-flat
and consider a short exact sequence of $kG$-modules
\[ 0 \longrightarrow K \longrightarrow L \longrightarrow M
     \longrightarrow 0 , \]
as in Theorem 3.1(iii). The double-character $kG$-module
$D^2K$ being cotorsion, the canonical map
$\jmath : K \longrightarrow D^2K$ factors through $L$. Since
$\jmath$ is a pure monomorphism (see, for example,
\cite[$\S $II.1.1.5]{RG}), it follows readily that the
monomorphism $K \longrightarrow L$ is pure as well. Therefore,
$M$ is a pure quotient of the cofibrant module $L$. \hfill $\Box$

\vspace{0.1in}

\noindent
The following result shows that the equality
${\tt Cof.flat}(kG) = {\tt GFlat}(kG)$ is a weak version of
the equality ${\tt Cof}(kG) = {\tt GProj}(kG)$.

\begin{Proposition}
Consider the following conditions:

(i) All Gorenstein projective $kG$-modules are cofibrant.

(ii) All projectively coresolved Gorenstein flat $kG$-modules
are cofibrant.

(iii) All Gorenstein flat $kG$-modules are cofibrant-flat.
\newline
Then, (i)$\rightarrow$(ii)$\leftrightarrow$(iii).
\end{Proposition}
\vspace{-0.05in}
\noindent
{\em Proof.}
It is clear that (i)$\rightarrow$(ii).

We now assume that (ii) holds and consider a Gorenstein flat
$kG$-module $M$. Invoking \cite[Theorem 4.11(4)]{SS}, we fix
a short exact sequence of $kG$-modules
\[ 0 \longrightarrow M \longrightarrow F \longrightarrow N
     \longrightarrow 0 , \]
where $F$ is flat and $N$ is projectively coresolved Gorenstein
flat. In view of assumption (ii), the $kG$-module $N$ is cofibrant
and hence cofibrant-flat. Since this is also the case for the flat
$kG$-module $F$, the closure of ${\tt Cof.flat}(kG)$ under
kernels of epimorphisms (cf.\ Proposition 2.2(ii)) implies
that the $kG$-module $M$ is cofibrant-flat. Hence, (iii) holds.

Conversely, if we assume that (iii) holds, then Corollary 2.6
implies that
\[ {\tt Cof}(kG) = {\tt PGF}(kG) \cap {\tt Cof.flat}(kG)
                 = {\tt PGF}(kG) \cap {\tt GFlat}(kG)
                 = {\tt PGF}(kG) , \]
so that (ii) holds. \hfill $\Box$

\vspace{0.1in}

\noindent
{\sc II.\ Two Hovey triples.}
The category ${\tt Cof}(kG)$ of cofibrant modules is a Frobenius
category with projective-injective objects the projective
$kG$-modules; cf.\ \cite[Lemma 4.5]{ER}. The following result is
a version of \cite[Theorem 19]{DE} for $kG$-modules.

\begin{Proposition}
There is a hereditary Hovey triple
$\left( {\tt Cof}(kG),{\tt Flat}(kG),{\tt Cof.flat}(kG) \right)$
in the exact category ${\tt Cof.flat}(kG)$. The homotopy category
of the associated exact model structure is equivalent, as a
triangulated category, to the stable category of cofibrant modules.
\end{Proposition}
\vspace{-0.05in}
\noindent
{\em Proof.}
The class ${\tt Flat}(kG)$ of flat $kG$-modules is closed under
direct summands and has the 2-out-of-3 property for short exact
sequences within the class of cofibrant-flat modules. Indeed,
${\tt Flat}(kG)$ is closed under extensions and kernels of
epimorphisms. On the other hand, if the cokernel of a monomorphism
between flat $kG$-modules is cofibrant-flat, then that cokernel is
necessarily flat. The latter claim follows, since any Gorenstein flat
$kG$-module of finite flat dimension is flat; cf.\ \cite[Corollary
10.3.4]{EJ}. We shall now prove that the pairs
\[ \left( {\tt Cof}(kG), {\tt Flat}(kG) \cap {\tt Cof.flat}(kG)
   \right) \;\; \mbox{and} \;\; \left(
   {\tt Cof}(kG) \cap {\tt Flat}(kG),{\tt Cof.flat}(kG) \right) \]
are complete and hereditary cotorsion pairs in ${\tt Cof.flat}(kG)$.
As we noted in Remark 3.4(ii), the intersection
${\tt Cof}(kG) \cap {\tt Flat}(kG)$ is equal to the class
${\tt Proj}(kG)$ of projective $kG$-modules; hence, the two
pairs displayed above become
\[ \left( {\tt Cof}(kG) , {\tt Flat}(kG) \right)
   \;\; \mbox{and} \;\;
   \left( {\tt Proj}(kG) , {\tt Cof.flat}(kG) \right) . \]
Corollary 3.3 implies that the pair
$\left( {\tt Cof}(kG),{\tt Flat}(kG) \right)$ is indeed a cotorsion
pair in the exact category ${\tt Cof.flat}(kG)$. The cotorsion pair
is complete, in view of the exact sequences obtained in Theorem
3.1(ii),(iii), and hereditary, since ${\tt Cof}(kG)$ is closed
under kernels of epimorphisms; cf.\ \cite[Proposition 2.2(ii)]{ER}.
Regarding the second pair diaplayed above, we note that
${\tt Cof.flat}(kG)$ is obviously the right Ext$^1$-orthogonal of
${\tt Proj}(kG)$ within ${\tt Cof.flat}(kG)$. In order to prove that
${\tt Proj}(kG)$ is the left Ext$^1$-orthogonal of ${\tt Cof.flat}(kG)$
within ${\tt Cof.flat}(kG)$, we consider a cofibrant-flat module $M$
which is contained in $^{\perp}{\tt Cof.flat}(kG)$ and fix a short
exact sequence of $kG$-modules
\[ 0 \longrightarrow M' \longrightarrow P \longrightarrow M
     \longrightarrow 0 , \]
where $P$ is projective. Since $M'$ is also cofibrant-flat, in view
of Proposition 2.2(ii), that short exact sequence splits and hence
the $kG$-module $M$ is projective. The closure of ${\tt Cof.flat}(kG)$
under kernels of epimorphisms implies the completeness of the cotorsion
pair. Of course, the cotorsion pair is hereditary.

The final claim in the statement of the Proposition follows from
\cite[Proposition 4.4 and Corollary 4.8]{G}. \hfill $\Box$

\vspace{0.1in}

\noindent
The exact category ${\tt Cotor}(kG) \cap {\tt Cof.flat}(kG)$
of cotorsion cofibrant-flat modules is Frobenius with
projective-injective objects the flat cotorsion $kG$-modules.
Indeed, we know that all flat cotorsion $kG$-modules are
contained in ${\tt Cof.flat}(kG)^{\perp}$ and hence these are
injective objects in ${\tt Cotor}(kG) \cap {\tt Cof.flat}(kG)$.
Therefore, Proposition 2.2(iv) implies that the exact category
${\tt Cotor}(kG) \cap {\tt Cof.flat}(kG)$ has enough injective
objects and all of these objects are necessarily flat cotorsion
$kG$-modules. On the other hand, the flat cotorsion $kG$-modules
are obviously projective objects in
${\tt Cotor}(kG) \cap {\tt Cof.flat}(kG)$. Using again Proposition
2.2(iv), it follows that ${\tt Cotor}(kG) \cap {\tt Cof.flat}(kG)$
has enough projective objects and all of these objects are flat
cotorsion $kG$-modules. The following result shows that we may
realize the stable category of
${\tt Cotor}(kG) \cap {\tt Cof.flat}(kG)$ as the homotopy category
of the exact model structure induced by a hereditary Hovey triple
in the category of $kG$-modules. That Hovey triple is a version of
the Hovey triple described by \v{S}aroch and
\v{S}t$\!$'$\!$ov\'{i}\v{c}ek, following \cite[Corollary 4.12]{SS},
for $kG$-modules.

\begin{Proposition}
There is a hereditary Hovey triple
$\left( {\tt Cof.flat}(kG),{\tt Cof}(kG)^{\perp},{\tt Cotor}(kG)
 \right)$
in the category of $kG$-modules. The homotopy category of the
associated model structure is equivalent to the stable category
of the Frobenius category ${\tt Cotor}(kG) \cap {\tt Cof.flat}(kG)$.
\end{Proposition}
\vspace{-0.05in}
\noindent
{\em Proof.}
The class ${\tt Cof}(kG)^{\perp}$ is closed under direct summands
and has the 2-out-of-3 property for short exact sequences; cf.\
\cite[Proposition 3.4]{ER}. We have to prove that the pairs
\[ \left( {\tt Cof.flat}(kG),{\tt Cof}(kG)^{\perp} \cap{\tt Cotor}(kG)
 \right) \;\; \mbox{and} \;\; \left(
 {\tt Cof.flat}(kG) \cap {\tt Cof}(kG)^{\perp},{\tt Cotor}(kG) \right) \]
are complete and hereditary cotorsion pairs in the category of $kG$-modules.
Indeed, Corollary 3.2 implies that the first one of these pairs is precisely
the cotorsion pair of Theorem 2.8. On the other hand, Corollary 3.3(ii)
shows that the second pair displayed above is precisely the flat cotorsion
pair $\left( {\tt Flat}(kG),{\tt Cotor}(kG) \right)$. The final statement
follows from \cite[Proposition 4.4 and Corollary 4.8]{G}. \hfill $\Box$

\section{Fibrant modules}

\noindent
In this section, we consider the class of fibrant modules and
obtain some basic properties of that class. We relate fibrant
modules to Gorenstein injective modules and show that they
form the right hand side of a complete hereditary cotorsion
pair. We also explore the duality between cofibrant-flat and
fibrant modules.

\vspace{0.1in}

\noindent
{\sc I.\ Fibrant and Gorenstein injective modules.}
We say that a $kG$-module $M$ is fibrant if the (diagonal)
$kG$-module $\mbox{Hom}_k(B,M)$ is injective and denote by
${\tt Fib}(kG)$ the class of fibrant modules. As any injective
$kG$-module is fibrant, the following conditions are equivalent
for a $kG$-module $N$:

(i) There is a non-negative integer $n$ and an exact sequence
of $kG$-modules
\[ 0 \longrightarrow N \longrightarrow M^0
     \longrightarrow M^1 \longrightarrow \cdots
     \longrightarrow M^n \longrightarrow 0 , \]
such that $M^0,M^1, \ldots ,M^n \in {\tt Fib}(kG)$.

(ii) $\mbox{id}_{kG}\mbox{Hom}_k(B,N) < \infty$.
\newline
If these two conditions are satisfied, we say that $N$ has finite
fibrant dimension; we denote by $\overline{\tt Fib}(kG)$ the class
of these $kG$-modules. It is clear that
${\tt Fib}(kG) \subseteq \overline{\tt Fib}(kG)$.

\begin{Lemma}\label{lem:Ben.inj}
Let $M \in \overline{\tt Fib}(kG)$ and assume that
${\rm Ext}^1_{kG}(N,M)=0$ for any $kG$-module $N$
of finite injective dimension. Then, $M \in {\tt Fib}(kG)$.
\end{Lemma}
\vspace{-0.05in}
\noindent
{\em Proof.}
Assume that $M$ is not fibrant, so that
$n = \mbox{id}_{kG}\mbox{Hom}_k(B,M)$ is strictly positive.
We consider a short exact sequence of $kG$-modules
\[ 0 \longrightarrow \mbox{Hom}_k(B,M)
     \stackrel{\imath}{\longrightarrow} I
     \longrightarrow J \longrightarrow 0 , \]
where $I$ is injective; of course, $J$ has finite injective
dimension. We also consider the commutative diagram with exact
rows
\[
\begin{array}{ccccccccc}
 0 & \longrightarrow
   & \mbox{Hom}_k \! \left( \overline{B},M \right)
   & \longrightarrow & \mbox{Hom}_k(B,M) & \longrightarrow
   & M & \longrightarrow & 0 \\
 & & \parallel & & \!\!\! {\scriptstyle{\imath}} \downarrow
 & & \!\!\! {\scriptstyle{\jmath}} \downarrow & & \\
 0 & \longrightarrow
   & \mbox{Hom}_k \! \left( \overline{B},M \right)
   & \longrightarrow & I &\longrightarrow & L
   & \longrightarrow & 0
\end{array}
\]
where the top row is obtained by applying the functor
$\mbox{Hom}_k(\_\!\_,M)$ to the $k$-split short exact
sequence of $kG$-modules
\[  0 \longrightarrow k \longrightarrow B
      \longrightarrow \overline{B} \longrightarrow 0 . \]
Since $\imath$ is injective, the snake lemma implies that
$\jmath$ is also injective and
$\mbox{coker} \, \jmath = \mbox{coker} \, \imath = J$. Then,
our assumption on $M$ implies that the monomorphism $\jmath$
splits, so that $M$ is a direct summand of $L$. It follows
that $\mbox{Hom}_k(B,M)$ is a direct summand of $\mbox{Hom}_k(B,L)$.
This is a contradiction, since $\mbox{Hom}_k(B,M)$ has injective
dimension equal to $n$, whereas the short exact sequence
\[ 0 \longrightarrow \mbox{Hom}_k \! \left( B ,
     \mbox{Hom}_k \! \left( \overline{B},M \right) \right)
     \longrightarrow \mbox{Hom}_k(B,I)
     \longrightarrow \mbox{Hom}_k(B,L) \longrightarrow 0 \]
and the fact that
$\mbox{Hom}_k \! \left( B ,
 \mbox{Hom}_k \! \left( \overline{B},M \right) \right) \!
 \simeq
 \mbox{Hom}_k \! \left( B \otimes_k \overline{B},M \right) \!
 \simeq
 \mbox{Hom}_k \! \left( \overline{B}, \mbox{Hom}_k(B,M) \right)$
has injective dimension $\leq n$ imply that $\mbox{Hom}_k(B,L)$
has injective dimension $\leq n-1$. \hfill $\Box$

\vspace{0.1in}

\noindent
Fibrant modules are closely related to Gorenstein injective modules.
Recall that an acyclic complex of injective $kG$-modules is called
totally acyclic if it remains acyclic after applying the functor
$\mbox{Hom}_{kG}(I,\_\!\_)$ for any injective $kG$-module $I$. A
$kG$-module $M$ is called Gorenstein injective if it is a kernel
of a totally acyclic complex of injective $kG$-modules. Let
${\tt GInj}(kG)$ denote the class of Gorenstein injective
$kG$-modules. We refer to \cite{EJ, Hol} for this notion (which
may be developed over any associative ring). Using the dual version
of an elegant construction by Cornick and Kropholler \cite{CK}, it
is shown in \cite[Proposition 5.6(i)]{ET3} that any fibrant module
is Gorenstein injective\footnote{The proof of \cite[Proposition
5.6(i)]{ET3} is given in the special case where the coefficient
ring is $\mathbb{Z}$. The same argument though applies verbatim
over any coefficient ring $k$.}, so that
${\tt Fib}(kG) \subseteq {\tt GInj}(kG)$. In particular,
\cite[Theorem 2.22]{Hol} implies that the functor
${\rm Ext}^1_{kG}(N,\_\!\_)$ vanishes on fibrant modules, whenever
$N$ is a $kG$-module of finite injective dimension. Hence, the
sufficient condition described in Lemma 4.1 for a $kG$-module to be
fibrant is also necessary. We list some basic properties of the class
${\tt Fib}(kG)$ of fibrant modules.

\begin{Proposition}\label{prop:fib}
(i) ${\tt Fib}(kG) = \overline{\tt Fib}(kG) \cap {\tt GInj}(kG)$.

(ii) The class ${\tt Fib}(kG)$ is closed under extensions, cokernels
of monomorphisms, direct products and direct summands.

(iii) Any fibrant module is a kernel of an acyclic complex of injective
$kG$-modules, all of whose kernels are fibrant.
\end{Proposition}
\vspace{-0.05in}
\noindent
{\em Proof.}
(i) We know that ${\tt Fib}(kG) \subseteq \overline{\tt Fib}(kG)$
and ${\tt Fib}(kG) \subseteq {\tt GInj}(kG)$. On the other hand,
if $M \in \overline{\tt Fib}(kG) \cap {\tt GInj}(kG)$, then
${\rm Ext}^1_{kG}(N,M)=0$ for any $kG$-module $N$ of finite
injective dimension; cf.\ \cite[Theorem 2.22]{Hol}. Hence,
Lemma 4.1 implies that $M$ is fibrant.

(ii) This follows since the class of injective $kG$-modules is
closed under extensions, cokernels of monomorphisms, direct
products and direct summands.

(iii) Let $N$ be a fibrant module; in particular, $N$ is Gorenstein
injective. Then, there exists a totally acyclic complex $I$ of
injective $kG$-modules, such that $N = Z^0I$. Since all injective
$kG$-modules are fibrant, it follows from (ii) above that the kernels
$(Z^nI)_n$ are fibrant for all $n \geq 0$. If $n<0$, then $Z^nI$
is a Gorenstein injective $kG$-module of finite fibrant dimension;
in view of (i) above, we conclude that $Z^nI \in {\tt Fib}(kG)$ as
well. \hfill $\Box$

\vspace{0.1in}

\noindent
Let us denote by $\mathcal{I}(kG)$ the class of those $kG$-modules
that appear as kernels of acyclic complexes of injective $kG$-modules.
Then, ${\tt GInj}(kG) \subseteq \mathcal{I}(kG)$ and equality holds
if and only if any acyclic complex of injective $kG$-modules is
totally acyclic. As a dual to the conjecture by Dembegioti and Talelli
\cite{DT} about cofibrant modules, we may ask whether
$\mathcal{I}(kG) \subseteq {\tt Fib}(kG)$, in which case we have
equalities ${\tt Fib}(kG) = {\tt GInj}(kG) = \mathcal{I}(kG)$.
We recall from $\S $1.II that $\mathfrak{Y}$ is the smallest class
of groups that contains all groups of type $\Phi$ over $k$ and
satisfies the equalities
${\scriptstyle{{\bf LH}}}\mathfrak{Y} = \mathfrak{Y}$ and
$\Phi_{inj} \mathfrak{Y} = \mathfrak{Y}$. The class $\mathfrak{Y}$
contains all ${\scriptstyle{{\bf LH}}}\mathfrak{F}$-groups.

\begin{Proposition}\label{prop:F=GInj}
Assume that any acyclic complex of injective $k$-modules is
contractible\footnote{This is the case if, for example, $k$
has finite weak global dimension; cf.\ \cite[Corollary 5.9]{Sto}.}
and let $G$ be a $\mathfrak{Y}$-group. Then,
${\tt Fib}(kG) = {\tt GInj}(kG) = \mathcal{I}(kG)$.
\end{Proposition}
\vspace{-0.05in}
\noindent
{\em Proof.}
We shall prove that $\mathcal{I}(kG) \subseteq {\tt Fib}(kG)$.
As in the proof of Proposition 2.3, we denote by $\mathcal{B}(kG)$
the class consisting of those $kG$-modules $B'$ which are projective
over any finite subgroup of $G$ and let
${\tt Fib}'(kG)$ be the class consisting of those $kG$-modules
$M$, for which the diagonal $kG$-module $\mbox{Hom}_k(B',M)$ is
injective for all $B' \in \mathcal{B}(kG)$. Since
$B \in \mathcal{B}(kG)$, it is clear that
${\tt Fib}'(kG) \subseteq {\tt Fib}(kG)$. Hence, it suffices to
show that $\mathcal{I}(kG) \subseteq {\tt Fib}'(kG)$. To that
end, we consider the class
\[ \mathfrak{W} =
   \{ H : \mathcal{I}(kH) \subseteq {\tt Fib}'(kH) \} \]
and show that $\mathfrak{Y} \subseteq \mathfrak{W}$. In view of
the definition of $\mathfrak{Y}$, it suffices to show that:
\newline
(i) all groups of type $\Phi$ over $k$ are contained in
$\mathfrak{W}$ and
\newline
(ii) ${\scriptstyle{{\bf LH}}}\mathfrak{W} = \mathfrak{W}$ and
$\Phi_{inj}\mathfrak{W} = \mathfrak{W}$.

In order to prove (i), we consider a group $H$ of type $\Phi$ over
$k$ and let $M \in \mathcal{I}(kH)$. Since injective $kH$-modules
are also $k$-injective, our assumption on $k$ implies that $M$ is
$k$-injective. We also consider a $kH$-module $B' \in \mathcal{B}(kH)$
and note that $B'$ has finite projective dimension, say
$\mbox{pd}_{kH}B' = n$. It follows that
$\mbox{id}_{kH}\mbox{Hom}_{kH}(B',M) \leq \mbox{pd}_{kH}B'=n$.
The latter inequality is valid for any $M \in \mathcal{I}(kH)$.
Since $M$ is the $n$-th cosyzygy of another $kH$-module
$M' \in \mathcal{I}(kH)$ and $B'$ is $k$-flat (being $k$-projective),
we may conclude that the $kH$-module $\mbox{Hom}_{kH}(B',M)$ is
actually injective. It follows that $M \in {\tt Fib}'(kH)$ and
hence $H \in \mathfrak{W}$.

We shall now prove (ii). In order to verify that $\mathfrak{W}$
is ${\scriptstyle{{\bf LH}}}$-closed, we consider an
${\scriptstyle{{\bf LH}}}\mathfrak{W}$-group $H$ and let
$M \in \mathcal{I}(kH)$ and $B' \in \mathcal{B}(kH)$. Then, we can
show that the $kH$-module $\mbox{Hom}_k(B',M)$ is injective, following
verbatim the arguments in the proof of \cite[Proposition 5.6(ii)]{ET3}
and proving that $\mbox{Hom}_k(B',M)$ is injective over any
${\scriptstyle{{\bf LH}}}\mathfrak{W}$-subgroup of $H$. The first step
of the transfinite induction therein (showing that $\mbox{Hom}_k(B',M)$
is injective over any ${\scriptstyle{{\bf H}}}\mathfrak{W}$-subgroup of
$H$) follows from the definition of $\mathfrak{W}$. We conclude that
$M \in {\tt Fib}'(kH)$, so that $H \in \mathfrak{W}$.

Finally, we shall prove that $\mathfrak{W}$ is $\Phi_{inj}$-closed.
To that end, we consider a $\Phi_{inj}\mathfrak{W}$-group $H$ and
let $M \in \mathcal{I}(kH)$ and $B' \in \mathcal{B}(kH)$. Then, the
definition of $\mathfrak{W}$ implies that $\mbox{Hom}_k(B',M)$ is
injective over any $\mathfrak{W}$-subgroup of $H$. Since $H$ is a
$\Phi_{inj}\mathfrak{W}$-group, it follows that the $kH$-module
$\mbox{Hom}_k(B',M)$ has finite injective dimension. This is the
case for any $M \in \mathcal{I}(kH)$; in other words, having fixed
$B' \in \mathcal{B}(kH)$, the injective dimension of
$\mbox{Hom}_k(B',M')$ is finite for any $M' \in \mathcal{I}(kH)$.
Since the class $\mathcal{I}(kH)$ is closed under direct products,
it is easily seen that there is an upper bound, say $n$, for the
injective dimensions $\mbox{id}_{kH}\mbox{Hom}_k(B',M')$,
$M' \in \mathcal{I}(kH)$. We note that $M$ is an $n$-th cosyzygy
of another $kG$-module $M' \in \mathcal{I}(kG)$ and $B'$ is $k$-flat;
we thus conclude that the $kG$-module $\mbox{Hom}_k(B',M)$ is actually
injective. It follows that $M \in {\tt Fib}'(kH)$ and hence
$H \in \mathfrak{W}$. \hfill $\Box$

\vspace{0.1in}

\noindent
{\bf Remarks 4.4.}
(i) Assuming the equality ${\tt Fib}(kG) = {\tt GInj}(kG)$, Lemma
\ref{lem:Ben.inj} states that, within the class of $kG$-modules
of finite Gorenstein injective dimension, the Gorenstein injective
$kG$-modules are those $kG$-modules which are right
Ext$^1$-orthogonal to the $kG$-modules of finite injective
dimension; see \cite[Theorem 2.15]{Hol}. In the same way, the
equality ${\tt Fib}(kG) = {\tt GInj}(kG)$ reduces the assertions
in Proposition \ref{prop:fib}(ii) to well-known properties of
Gorenstein injective modules (cf.\ \cite[Theorem 2.6]{Hol}),
whereas Proposition \ref{prop:fib}(iii) is a consequence of the
very definition of Gorenstein injective modules.

(ii) For any $kG$-module $M$ we define its fibrant dimension
as the injective dimension of the $kG$-module $\mbox{Hom}_k(B,M)$.
Since fibrant modules are Gorenstein injective, this dimension
bounds the Gorenstein injective dimension of $M$. The global
fibrant dimension $\mbox{gl.fib.dim} \, kG$ of $kG$ is then
defined as the supremum of the fibrant dimensions of all
$kG$-modules. Analogously, for any $kG$-module $N$ we define
its cofibrant dimension as the projective dimension of the
$kG$-module $B \otimes_kN$. Since cofibrant modules are
Gorenstein projective, this dimension bounds the Gorenstein
projective dimension of $N$. Then, the global cofibrant dimension
$\mbox{gl.cof.dim} \, kG$ of $kG$ is the supremum of the
cofibrant dimensions of all $kG$-modules. The standard
Hom-tensor adjunction for diagonal $kG$-modules implies that
there is an isomorphism of abelian groups
\[ {\rm Ext}_{kG}^n(B \otimes_kN,M) \simeq
   {\rm Ext}_{kG}^n \! \left( N, \mbox{Hom}_k(B,M) \right) \]
for any $kG$-modules $N,M$ and any $n \geq 0$. It follows that
$\mbox{gl.cof.dim} \, kG = \mbox{gl.fib.dim} \, kG$.
\addtocounter{Lemma}{1}

\vspace{0.1in}

\noindent
{\sc II.\ The cotorsion pair
$(^{\perp}{\tt Fib}(kG),{\tt Fib}(kG))$.}
Baer's criterion for injectivity implies that a $kG$-module is
injective if and only if it is right Ext$^1$-orthogonal to all
cyclic $kG$-modules. Using this fact, we can prove the following
result.

\begin{Theorem}\label{thm:cp-inj}
The pair $\left( ^{\perp}{\tt Fib}(kG),{\tt Fib}(kG) \right)$
is a hereditary cotorsion pair in the category of $kG$-modules,
which is cogenerated by a set. Its kernel is the class of
injective $kG$-modules.
\end{Theorem}
\vspace{-0.05in}
\noindent
{\em Proof.}
Let $\mathcal{T}$ be a {\em set} of representatives of the
isomorphism classes of all cyclic $kG$-modules and consider
$\mathcal{S} = \{ C \otimes_kB : C \in \mathcal{T} \}$. If $M$
is a $kG$-module and $C \in \mathcal{T}$, then the isomorphism
${\rm Ext}_{kG}^1(C \otimes_kB,M) \simeq
 {\rm Ext}_{kG}^1 \! \left( C, \mbox{Hom}_k(B,M) \right)$
and Baer's injectivity criterion imply that $M$ is fibrant if
and only if $M \in \mathcal{S}^{\perp}$. Therefore,
${\tt Fib}(kG) = \mathcal{S}^{\perp}$ and hence
\[ \left( ^{\perp}{\tt Fib}(kG),{\tt Fib}(kG) \right) =
   \left( ^{\perp} \! \left( \mathcal{S}^{\perp} \right) \! ,
   \mathcal{S}^{\perp} \right) \]
is the cotorsion pair cogenerated by the set $\mathcal{S}$.
The cotorsion pair is hereditary, since the right hand class
${\tt Fib}(kG)$ is closed under cokernels of monomorphisms;
cf.\ Proposition 4.2(ii). We shall prove that the kernel of
the cotorsion pair is the class of injective $kG$-modules,
i.e.\ that
\[ ^{\perp}{\tt Fib}(kG) \cap {\tt Fib}(kG) = {\tt Inj}(kG) . \]
Let $M$ be a $kG$-module contained in the intersection
$^{\perp}{\tt Fib}(kG) \cap {\tt Fib}(kG)$. Then, $M$ is fibrant and
hence Proposition 4.2(iii) implies that there exists a short exact
sequence of $kG$-modules
\[  0 \longrightarrow N \longrightarrow I \longrightarrow M
      \longrightarrow 0 , \]
where $N$ is fibrant and $I$ is injective. Since
$M \in \, \! ^{\perp}{\tt Fib}(kG)$, that exact sequence splits
and hence $M$ is injective. Conversely, we note that all injective
$kG$-modules are fibrant, whereas the inclusion
${\tt Fib}(kG) \subseteq {\tt GInj}(kG)$ implies that
${\tt Inj}(kG) \subseteq \, \! ^{\perp}{\tt GInj}(kG)
               \subseteq \, \! ^{\perp}{\tt Fib}(kG)$.
\hfill $\Box$

\vspace{0.1in}

\noindent
We now list a few properties of the class $^{\perp}{\tt Fib}(kG)$.

\begin{Proposition}\label{prop:fib+}
The class $^{\perp}{\tt Fib}(kG)$ is closed under direct summands
and has the 2-out-of-3 property for short exact sequences of
$kG$-modules. It contains all $kG$-modules of finite flat or
injective dimension and all $kG$-modules of the form
$M \otimes_kB$, where $M$ is a $kG$-module.
\end{Proposition}
\vspace{-0.05in}
\noindent
{\em Proof.}
The class $^{\perp}{\tt Fib}(kG)$ is clearly closed under direct
summands and extensions. Since the cotorsion pair
$\left( ^{\perp}{\tt Fib}(kG) , {\tt Fib}(kG) \right)$
is hereditary, the class $^{\perp}{\tt Fib}(kG)$ is closed
under kernels of epimorphisms. Consider a short
exact sequence of $kG$-modules
\begin{equation}
 0 \longrightarrow M' \longrightarrow M \longrightarrow M''
   \longrightarrow 0 ,
\end{equation}
where $M',M \in \, \! ^{\perp}{\tt Fib}(kG)$. We shall prove
that $M'' \in \, \! ^{\perp}{\tt Fib}(kG)$. We fix a fibrant
module $N$ and invoke Proposition \ref{prop:fib}(iii), in order
to find a short exact sequence of $kG$-modules
\begin{equation}
 0 \longrightarrow N' \longrightarrow I \longrightarrow N
   \longrightarrow 0 ,
\end{equation}
where $N'$ is fibrant and $I$ is injective. Since both groups
${\rm Ext}^1_{kG}(M',N')$ and ${\rm Ext}^2_{kG}(M,N')$ are
trivial, it follows from the short exact sequence (5) that
${\rm Ext}^2_{kG}(M'',N')=0$. The short exact sequence (6)
implies then that
${\rm Ext}^1_{kG}(M'',N) = {\rm Ext}^2_{kG}(M'',N')$ is the
trivial group and hence $M'' \in \, \! ^{\perp}{\tt Fib}(kG)$,
as needed.

Since any fibrant module is Gorenstein injective, the class
$^{\perp}{\tt Fib}(kG)$ contains all $kG$-modules of finite
injective dimension (cf.\ \cite[Theorem 2.22]{Hol}) and all
$kG$-modules of finite flat dimension (cf.\ \cite[Corollary 5.9]{Sto}).
Finally, if $M$ is any $kG$-module, then the Hom-tensor
adjunction for diagonal $kG$-modules implies that the functor
\[ {\rm Ext}^1_{kG}(M \otimes_kB, \_\!\_) =
   {\rm Ext}^1_{kG}(M,\mbox{Hom}_k(B,\_\!\_)) \]
vanishes on all fibrant modules and hence
$M \otimes_kB \in \, \! ^{\perp}{\tt Fib}(kG)$. \hfill $\Box$

\vspace{0.1in}

\noindent
The exact category ${\tt Fib}(kG)$ of fibrant modules is
Frobenius with projective-injective objects the injective
$kG$-modules. Indeed, it is clear that all injective
$kG$-modules are injective objects in ${\tt Fib}(kG)$.
Proposition 4.2(iii) implies then that ${\tt Fib}(kG)$ has
enough injective objects and all of these objects are
injective $kG$-modules. On the other hand, all injective
$kG$-modules are projective objects in ${\tt Fib}(kG)$,
since fibrant modules are Gorenstein injective. Using again
Proposition 4.2(iii), we conclude that ${\tt Fib}(kG)$ has
enough projective objects and all of these objects are injective
$kG$-modules. Theorem 4.5 implies that there is a hereditary
Hovey triple
\begin{equation}
 \left( kG\mbox{-Mod}, \, \! ^{\perp}{\tt Fib}(kG) ,
 {\tt Fib}(kG) \right)
\end{equation}
in the category of $kG$-modules. The homotopy category of the
associated model structure is equivalent to the stable category
of fibrant modules.

\vspace{0.1in}

\noindent
{\bf Remarks 4.7.}
(i) Let $M$ be a $kG$-module of finite fibrant dimension. In
particular, $M$ has finite Gorenstein injective dimension and
hence \cite[Theorem 2.15]{Hol} implies that there exists a short
exact sequence of $kG$-modules
\[ 0 \longrightarrow M \longrightarrow J \longrightarrow N
     \longrightarrow 0 , \]
where $J$ is Gorenstein injective and $\mbox{id}_{kG}N < \infty$.
Since both $kG$-modules $M,N$ have finite fibrant dimension, the
same is true for $J$. Invoking Proposition 4.2(i), we conclude
that the $kG$-module $J$ is actually fibrant. Using standard
techniques, the existence of the above short exact sequence
implies the existence of another short exact sequence of $kG$-modules
\[ 0 \longrightarrow J' \longrightarrow N' \longrightarrow M
     \longrightarrow 0 , \]
where $J'$ is fibrant and $\mbox{id}_{kG}N' < \infty$.\footnote{We
may consider a short exact sequence
$0 \longrightarrow J' \longrightarrow I \longrightarrow J
   \longrightarrow 0$,
where $J'$ is fibrant and $I$ is injective (cf.\ Proposition
4.2(iii)) and form its pullback along the monomorphism
$M \longrightarrow J$.}

(ii) Denoting by $\overline{\tt Inj}(kG)$ the class of all
$kG$-modules of finite injective dimension, we note that the
approximation sequences obtained in (i) above imply that
$\overline{\tt Fib}(kG) \cap \overline{\tt Inj}(kG)^{\perp} =
 {\tt Fib}(kG)$
and
$\overline{\tt Fib}(kG) \cap \, \! ^{\perp}{\tt Fib}(kG) =
 \overline{\tt Inj}(kG)$.
It follows readily that
$\overline{\tt Inj}(kG) \cap {\tt Fib}(kG) = {\tt Inj}(kG)$
and hence
\[ \left( \overline{\tt Fib}(kG) ,
   \overline{\tt Inj}(kG) , {\tt Fib}(kG) \right) \]
is a Hovey triple in the exact category $\overline{\tt Fib}(kG)$
of $kG$-modules of finite fibrant dimension. The latter Hovey
triple can be also obtained by restricting to
$\overline{\tt Fib}(kG)$ the Hovey triple displayed in (7) above.
\addtocounter{Lemma}{1}

\vspace{0.1in}

\noindent
The standard duality between flat and injective $kG$-modules,
which is induced by the character (Pontryagin duality) functor
\cite{L}, extends to a duality between the classes of cofibrant-flat
and fibrant modules.

\begin{Proposition}
A $kG$-module is cofibrant-flat if and only if its character
module is fibrant.
\end{Proposition}
\vspace{-0.05in}
\noindent
{\em Proof.}
Let $M$ be a $kG$-module. We have to show that the $kG$-module
$M \otimes_kB$ is flat if and only if the $kG$-module
$\mbox{Hom}_k(B,DM)$ is injective. Since there is an isomorphism
of $kG$-modules $\mbox{Hom}_k(B,DM) \simeq D(M \otimes_kB)$, the
result follows from \cite{L}. \hfill $\Box$

\begin{Corollary}
There exists a $kG$-module $T$, such that the class
${\tt Cof.flat}(kG)$ consists precisely of those
$kG$-modules $M$ for which ${\rm Tor}^{kG}_1(T,M)=0$.
\end{Corollary}
\vspace{-0.05in}
\noindent
{\em Proof.}
Theorem 4.5 implies that there exists a $kG$-module $T$, for
which ${\tt Fib}(kG) = \{ T \} ^{\perp}$. Invoking Proposition
4.8, it follows that ${\tt Cof.flat}(kG)$ consists precisely
of those $kG$-modules $M$, for which $DM \in \{ T \}^{\perp}$.
In order to complete the proof, we note that the abelian group
$D\mbox{Tor}^{kG}_1(T,M) = {\rm Ext}^1_{kG}(T,DM)$ is trivial
if and only if $\mbox{Tor}^{kG}_1(T,M) = 0$. \hfill $\Box$

\vspace{0.1in}

\noindent
{\bf Remarks 4.10.}
(i) The proof of Theorem 4.5 shows that a $kG$-module $T$ with the
property stated in Corollary 4.9 is the (diagonal) $kG$-module
$S \otimes_kB$, where $S$ is the direct sum of a set of
representatives of the isomorphism classes of cyclic $kG$-modules.

(ii) The class of Gorenstein flat modules over a right coherent
ring $R$ admits a description as the class of roots of the functor
$\mbox{Tor}_1^R(T,\_\!\_)$, for a suitable right $R$-module $T$;
cf.\ \cite[Corollary 5.7]{SS}. Combining Proposition 2.3 with
Corollary 4.9, we obtain another class of rings (which are often
non-coherent) over which the same description of Gorenstein flat
modules is possible.
\addtocounter{Lemma}{1}

\vspace{0.1in}

\noindent
We conclude with an addendum to Proposition 3.6.

\begin{Corollary}
If all Gorenstein injective $kG$-modules are fibrant, then
all Gorenstein flat $kG$-modules are cofibrant-flat.
\end{Corollary}
\vspace{-0.05in}
\noindent
{\em Proof.}
If $M$ is a Gorenstein flat $kG$-module, then its character
module is Gorenstein injective; cf.\ \cite[Theorem 3.6]{Hol}.
In view of our assumption, that character module is fibrant.
Proposition 4.8 then implies that $M$ is cofibrant-flat, as
needed. \hfill $\Box$

\vspace{0.1in}

\noindent
{\em Acknowledgments.}
Wei Ren was supported by the National Natural Science Foundation
of China (No. 11871125).

\vspace{0.1in}

{\footnotesize \noindent Ioannis Emmanouil\\
Department of Mathematics, University of Athens, Athens 15784, Greece \\
E-mail: {\tt emmanoui$\symbol{64}$math.uoa.gr}}

\vspace{0.05in}

{\footnotesize \noindent Wei Ren\\
 School of Mathematical Sciences, Chongqing Normal University, Chongqing 401331, PR China\\
 E-mail: {\tt wren$\symbol{64}$cqnu.edu.cn}}

\end{document}